\documentclass[]{article}
\usepackage[margin=2.6cm]{geometry}
\usepackage{graphicx,amsfonts,subfigure}
\usepackage{amsmath,amssymb,setspace}
\usepackage{stmaryrd}
\usepackage{eufrak}
\usepackage{bbm}
\usepackage[mathscr]{euscript}
\usepackage{hyperref}
\usepackage{algorithmic}
\usepackage{leftidx}

\numberwithin{equation}{section}

\newtheorem{thm}{Theorem}[section]
\newtheorem{cor}[thm]{Corollary}

\newtheorem{lemma}[thm]{Lemma}

\newtheorem{theorem}[thm]{Theorem}
\newtheorem{remark}[thm]{Remark}
\newtheorem{proposition}[thm]{Proposition}
\newtheorem{definition}[thm]{Definition}
\newtheorem{corollary}[thm]{Corollary}

\newenvironment{proof}{{\bf Proof.}}{\hfill$\square$\vskip.5cm}

\newcommand{\R}{\mathbb{R}}

\newcommand{\lin}{{\rm KR}}

\title{Analysis and Approximation of the Canonical Polyadic Tensor Decomposition}

\author{Stefan Kindermann\footnote{Industrial Mathematics Institute, 
Johannes Kepler Universitat Linz, Altenbergerstrasse 69, A-4040 Linz, Austria,
 kindermann@indmath.uni-linz.ac.at} \, \, \, \, 
Carmeliza Navasca\footnote{Department of Mathematics, 
Clarkson University, Potsdam, NY, 13699, USA, cnavasca@clarkson.edu. 
C.N. is in part supported by National Science Foundation DMS-0915100.}  }

\date{\today}

\begin{document}

\maketitle
\begin{abstract}
\setcounter{section}{0}
We study the least-squares (LS) functional of the canonical polyadic (CP) tensor decomposition. 
Our approach is based on the elimination of one factor matrix which results in a reduced functional. 
The reduced functional is reformulated into a projection framework and into a Rayleigh quotient. An analysis of this functional leads to several conclusions: new sufficient conditions for the existence of minimizers of the LS functional, the existence of a critical point in the rank-one case, a heuristic explanation of "swamping" and computable bounds on the minimal value of the LS functional. The latter result leads to a simple algorithm
-- the Centroid Projection algorithm -- to compute suboptimal solutions of tensor decompositions. 
These suboptimal solutions are applied to iterative CP algorithms as initial guesses, yielding a method called centroid projection for canonical polyadic  (CPCP) decomposition which provides a significant speedup in our numerical experiments compared to the standard methods.

\vspace{8pt}
\vskip .2 cm
\end{abstract}

\section{Introduction}

In 1927, Hitchcock \cite{Hitch1,Hitch2} introduced the idea that a tensor is decomposable into a sum 
of a finite number of rank-one tensors.  Today, this decomposition is referred to as the canonical 
polyadic (CP) tensor decomposition (also known as CANDECOMP \cite{CarolChang} or PARAFAC \cite{Harshman}).
 CP tensor decomposition reduces a tensor to  a linear combination of rank-one tensors, i.e.
\begin{eqnarray}\label{CP}
(\mathcal{A})_{ijk}= \sum_{r=1}^R a_{ir} b_{jr} c_{kr}
\end{eqnarray}
where $\mathcal{A}\in \mathbb{R}^{I \times J \times K}$, 
$\bold{a_r} = (a_{ir})_{i=1}^I \in \mathbb{R}^{I}, \bold{b_r} = (b_{jr})_{j=1}^J
\in \mathbb{R}^{J}$ and $\bold{c_r} = (c_{kr})_{k=1}^K \in \mathbb{R}^{K}$. 
The column vectors $\bold{a_r}, \bold{b_r}$ and $\bold{c_r}$ form the so-called factor matrices 
$\bold{A} \in \mathbb{R}^{I \times R} $, $\bold{B} \in \mathbb{R}^{J \times R}$ and
 $\bold{C} \in \mathbb{R}^{K \times R}$.  
The tensorial rank \cite{Hitch2} is the minimum $R \in \mathbb{N}$ such that $\mathcal{T}$ can
 be expressed as a sum of $R$ rank-one tensors. 

The problem of interest is to find -- if it exists -- the \emph{best} approximate 
tensor representable in a CP format with a tensorial rank $R$ 
from a given (possibly noisy) tensor  $\mathcal{T} \in \mathbb{R}^{I \times J \times K}$. 
A standard approach for this task is  to minimize the norm of the residual tensor in the least-square sense:
\begin{eqnarray}\label{mainobj}
\mathfrak{J}(\bold{A},\bold{B},\bold{C}) = \frac{1}{2} \sum_{i,j,k} 
\left ((\mathcal{T})_{ijk} - \sum_{r=1}^R a_{ir} b_{jr} c_{kr} \right )^2 .
\end{eqnarray}

A popular iterative method for approximating the given tensor 
$\mathcal{T}$ via its factors $(\bold{A}, \bold{B}, \bold{C})$ is 
called the Alternating Least-Squares (ALS) technique. Independently, 
ALS was introduced by Carol and Chang \cite{CarolChang} and Harshman \cite{Harshman} in 1970. 
The ALS method is an application of the nonlinear block Gauss-Seidel algorithm \cite{LiKinNav} 
where the nonlinear optimization (\ref{mainobj}) is reduced into several
 least-squares subproblems which are solved iteratively with subsequent 
updates of the factor minimizer.  The ALS algorithm has been extensively applied to many problems across various engineering and science disciplines; e.g., see the survey paper \cite{Kolda} and the references therein. Despite the widespread popularity of ALS, it has its shortcomings. 
Problems can arise in degenerate problems and slow converging nondegenerate problems with CP solutions. To this end, methods like regularization techniques \cite{NaDeLatKind} and enhanced line search \cite{RaCo} are improvements of ALS. There are also several methods based on other techniques, such as, conjugate gradient \cite{PaateroEngine}
 and Schur decomposition \cite{StegDeLat} for CP decomposition; see the paper of Comon et al. \cite{Comon} on the survey of ALS methods.

In this paper, we analyze the minimization of the objective function (\ref{mainobj}) by eliminating one factor $\bold{A}$, reducing to a minimization over the factor minimizers $\bold{B}$ and $\bold{C}$, equivalent to the original one. 
Analysis of the reduced functional allows reformulations into several forms:  as a Rayleigh quotient type functional or as an weighted projection
onto the Khatri-Rao range of  $\bold{B}$ and $\bold{C}$. As a consequence, we prove sufficient conditions for the existence of the minimizer of (\ref{mainobj})
in terms of the rank of the Khatri-Rao matrix which substantiates well-known facts about the degeneracy case, like the diverging norms
of the factors and that the solution space is not closed  \cite{TheDutch,PP,SilLim}. Furthermore, for the special case of
rank-1 decomposition, we show -- using Morse theory -- the existence of a critical point which can lead to 
a halt of the ALS algorithms at nonextremal points. Poor convergence (\emph{swamping}) of the ALS algorithm can be attributed to the feasible set, the Khatri-Rao range of  $\bold{B}$ and $\bold{C}$, of the reduced objective functional.

Further analysis of the reduced objective functional provides upper and lower bounds. 
The minimizers of the upper bound turn out to be computable by linear algebra methods,
yielding an effective and simple algorithm (the Centroid Projection (CP) Algorithm) 
for computing suboptimal solutions  to (\ref{mainobj}). The suboptimal solutions may
serve as a initial guesses to any iterative CP-decomposition methods like ALS or other advanced 
algorithms. We will refer to this powerful combination as the CPCP method.
 In our numerical examples, the Centroid Projection  have
 shown to improve performance of several iterative methods for CP decomposition in comparison
 to the examples with random initial starters. Moreover, initialization of the upper 
bound minimizers  works well for CP decomposition with symmetries \cite{Stegeman}, that is, when at least
 two of the factors are identical.

\section{Preliminaries}\label{sec:prel}

We denote the scalars in $\mathbb{R}$ with lower-case letters $(a,b,\ldots)$ and the vectors
 with bold lower-case letters $(\bf{a},\bf{b},\ldots)$.  
The matrices are written as bold upper-case letters $(\bf{A}, \bf{B},\ldots)$ and the symbol 
for tensors are calligraphic letters $(\mathcal{A},\mathcal{B},\ldots)$. 
The subscripts represent the following scalars:  $\mathcal{(A)}_{ijk}=a_{ijk}$, 
$(\bold{A})_{ij}=a_{ij}$, $(\bold{a})_i=a_i$. The superscripts indicate the length 
of the vector or the size of the matrices. For example, $\bold{b}^{K}$ is a vector
 with length $K$ and $\bold{B}^{N \times K}$ is a $N \times K$ matrix. In addition, 
the lower-case superscripts on a matrix indicate the mode in which has been matricized. 

The order of a tensor refers to the cardinality of the index set.  
A matrix is a second-order tensor and a vector is a first-order tensor. 
The scalar product of $\mathcal{T}$, $\mathcal{R}  \in \mathbb{R}^{I \times J \times K}$
is defined as  
\[ \langle \mathcal{T},\mathcal{L} \rangle = \sum_{ijk} (\mathcal{T})_{ijk} (\mathcal{L})_{ijk}.  \]
The Frobenius norm of $\mathcal{A} \in \mathbb{R}^{I \times J \times K}$ is defined as
 \[ \Vert \mathcal{A} \Vert_F^{2} = \sum_{i=1}^I \sum_{j=}^J \sum_{k=1}^K \vert a_{ijk} \vert^2 =\langle \mathcal{A}, \mathcal{A}  \rangle \]
which is a direct extension of the Frobenius norm of a matrix. Furthermore we denote
by $\cdot$ the usual matrix product. 
\begin{definition}
The Khatri-Rao product of $\bold{A} \in \mathbb{R}^{I \times R}$ 
and $\bold{B} \in \mathbb{R}^{J \times R}$ is defined as
\begin{eqnarray*}
\bold{A} \odot \bold{B}= [\bold{a_1} \otimes \bold{b_1} ~ \bold{a_2} \otimes \bold{b_2}~ \ldots  ~ \bold{a_R} \otimes \bold{b_R}] \in \mathbb{R}^{IJ \times R}
\end{eqnarray*}
when $\bold{A}=[\bold{a_1}~ \bold{a_2}~ \ldots ~\bold{a_R}]$ 
and $\bold{B}=[\bold{b_1}~ \bold{b_2}~ \ldots ~\bold{b_R}]$.
\end{definition}
Here, $\bold{a} \otimes \bold{b}$ denotes the Kronecker product of two vectors $\bold{a} \in \R^I$,
$\bold{b} \in \R^J$ yielding a vector of size $I J$ with entries that are all possible products 
of the entries in $\bold{a}$ and $\bold{b}$. 

\begin{definition}[Tucker mode-$n$ product]
Given a tensor $\mathcal{T} \in \mathbb{R}^{I_1 \times I_2 \times I_3}$ and matrices $\bold{A} \in \mathbb{R}^{I_1 \times J_1}$, $\bold{B} \in \mathbb{R}^{ I_2 \times J_2}$ and $\bold{C}\in \mathbb{R}^{I_3 \times J_3}$, then the Tucker mode-$n$ products are the following:
{\small \begin{eqnarray*}
\mathcal{T} \bullet_1 \bold{A} := (\mathcal{T} \bullet_1 \bold{A})_{j_1i_2i_3} &=& 
\sum_{i_1=1}^{I_1} t_{i_1i_2i_3}a_{i_1j_1},~\forall j_1,i_2,i_3~\hspace{.15cm}\mbox{(mode-1 product)}\\
\mathcal{T} \bullet_2 \bold{B} := (\mathcal{T} \bullet_2 \bold{B})_{i_1j_2i_3} &=& \sum_{i_2=1}^{I_2} t_{i_1i_2i_3}b_{i_2j_2},~\forall j_2,i_1,i_3~\hspace{.15cm}\mbox{(mode-2 product)}\\
\mathcal{T} \bullet_3 \bold{C} := (\mathcal{T} \bullet_3 \bold{C})_{i_1i_2j_3} &=& 
\sum_{i_3=1}^{I_3} t_{i_1i_2i_3}c_{i_3j_3},~\forall j_3,i_1,i_2~\hspace{.15cm}\mbox{(mode-3 product)}.
\end{eqnarray*}}
\end{definition}
Moreover, the Tucker mode products can be combined as in this example:
\[ \mathcal{T} \bullet_{2,3} (\bold{B}, \bold{C}):= (\mathcal{T} \bullet_{2,3} (\bold{B}, \bold{C}))_{i_1r}:= \sum_{i_2=1}^{I_2} \sum_{i_3=1}^{I_3} \mathcal{T}_{i_1i_2i_3} b_{i_2r} c_{i_3r}  \]
where $\bold{B} \in \mathbb{R}^{I_2 \times R}$ and $\bold{C} \in \mathbb{R}^{I_3 \times R}$.

\begin{definition}[outer product of vectors]
For vectors $\bold{a} \in \R^I$, $\bold{b} \in \R^J$ the outer product $\bold{a} \circ \bold{b}$
is the $I \times J$ matrix with entries
\[ (\bold{a} \circ \bold{b})_{i,j} = a_i b_j,~\forall i,j \]
similarly, the outer product of three vectors $\bold{a} \in \R^I$, $\bold{b} \in \R^J$, $\bold{c} \in \R^K$
is the $I \times J \times K$ tensor
\[ (\bold{a} \circ \bold{b} \circ \bold{c})_{i,j,k} =  a_i b_j c_k,~\forall i,j,k \]
\end{definition}
\section{The least squares functional and its reduction}
Recall the least-squares objective functional in (\ref{mainobj}):
\begin{eqnarray} \label{mainobj2}
\mathfrak{J}(\bold{A},\bold{B},\bold{C}) = \frac{1}{2}  
\left \Vert \mathcal{T} - \sum_{r=1}^R \bold{a_r} \circ \bold{b_r} \circ \bold{c_r}  \right \Vert _F^{2}
\end{eqnarray}
where $\Vert \cdot \Vert _F$ is the Frobenius norm. 
The goal is to find  minimizers $\bold{A}$, $\bold{B}$ and $\bold{C}$ of
\[ \inf_{\bold{A},\bold{B},\bold{C}} \mathfrak{J}(\bold{A},\bold{B},\bold{C}). \]
Note that it is well-known that this infimum is not necessarily attained see, e.g., \cite{SilLim}.

\begin{lemma}\label{fixedBC} 
Let $\bold{B},\bold{C}$ be fixed.  The solution to the minimization problem
\begin{equation}\label{defamin}  
\tilde{\bold{A}}[\bold{B},\bold{C}] := \mathrm{argmin}_{ \bold{A} \in \R^{I \times R}}
 \mathfrak{J}(\bold{A},\bold{B},\bold{C})  
\end{equation}
exists. In fact, a minimizer is given by
\begin{equation}\label{theA}
 \tilde{\bold{A}}[\bold{B},\bold{C}] =  \mathcal{T}\bullet_{2,3} (\bold{B},\bold{C}) \cdot \bold{G}^\dagger,
\end{equation}
where  $\bold{G}^\dagger$ is the pseudo-inverse of $\bold{G}$ with elements
\begin{eqnarray}\label{theG}
(\bold{G})_{rs} :=  \left (\sum_{j=1}^J b_{jr} b_{js} \right ) \left ( \sum_{k=1}^J c_{kr} c_{ks} \right ) 
\end{eqnarray}
and
\[ (\mathcal{T} \bullet_{2,3} (\bold{B}, \bold{C}))_{ir}:= \sum_{j=1}^J \sum_{k=1}^K \mathcal{T}_{ijk} b_{jr} c_{kr}.  \]
\end{lemma}
\begin{proof}
With 
$\bold{B},\bold{C}$ being fixed,  (\ref{defamin}) is a usual finite dimensional  linear least squares problem 
for which  it is well-known that a solution exists. 
Differentiation of the functional (\ref{mainobj}) with respect to $a_{i^*r^*}$ leads to
the optimality conditions
\[ \frac{\partial}{\partial a_{i^*r^*}} \mathfrak{J}(\bold{A},\bold{B},\bold{C}) = 
 \sum_{j,k} \left(\mathcal{T}_{i^*jk} - \sum_{r=1}^R a_{i^*r} b_{jr} c_{kr} \right) b_{jr^*} c_{kr^*} = 0.
\]
Since $\sum_{j,k} \left(\mathcal{T}_{i^*jk} - \sum_{r=1}^R a_{i^*r} b_{jr} c_{kr} \right) 
b_{jr^*} c_{kr^*}=
\sum_{j,k} \mathcal{T}_{i^*jk} b_{jr^*} c_{kr^*} - 
\sum_{j,k} \sum_{r=1}^R a_{i^*r} b_{jr} c_{kr} b_{jr^*} c_{jr^*}$, we obtain the matrix equation
\[ \bold{A}\cdot \bold{G} = \mathcal{T}\bullet_{2,3} (\bold{B},\bold{C}) \]
We know that a (not necessarily unique) solution exists, which then 
is expressible in terms of the pseudo-inverse (\ref{theA}).
\end{proof}
From the definition (\ref{theG}),
$$
(\bold{G})_{rs} = \left (\sum_{j=1}^J b_{jr} b_{js} \right ) \left ( \sum_{k=1}^K c_{kr} c_{ks} \right )
=(\bold{b_r} \otimes \bold{c_r})^T (\bold{b_s} \otimes \bold{c_s}),
$$
it follows that 
\begin{eqnarray}\label{theG2}
\bold{G}=(\bold{B} \odot \bold{C})^T \cdot (\bold{B} \odot \bold{C}) \in \mathbb{R}^{R \times R}, 
\end{eqnarray}
is a Gramian matrix for the vectors $\bold{b_r} \otimes \bold{c_r}$, $r=1,\ldots R$ as well as the 
Hadamard product of $\bold{B}^T\cdot \bold{B}$ and $\bold{C}^T\cdot \bold{C}$.  Note that 
$\bold{G}$ depends on $\bold{B}$ and $\bold{C}$ but we omitted 
this dependence to avoid exuberant notation. It follows easily that $\bold{G}$ is symmetric,
and thus $\bold{G}^\dagger$ is. Moreover, the pseudo-inverse satisfies the Moore-Penrose
equation $\bold{G}^{\dagger} \cdot \bold{G} \cdot \bold{G}^{\dagger}= \bold{G}^{\dagger}$.

%

Motivated by the ALS algorithm, which iteratively minimizes over the factors matrices,
we state the main tool in this paper, the reduced functional.
Minimization over $\bold{A}$ reduces the original functional so that  we now define
\begin{equation}\label{defJred}
\mathfrak{J}_{red}(\bold{B},\bold{C}):= \mathfrak{J}(\tilde{\bold{A}}[\bold{B},\bold{C}],\bold{B},\bold{C})
\end{equation} 
where $\tilde{\bold{A}}[\bold{B},\bold{C}]$ is a minimizer in \eqref{defamin}. 
This definition does
not depend which minimizer we take.
In the following lemma, we 
show that the minimizers of $\mathfrak{J}$ can be found through the minimizers of $\mathfrak{J}_{red}$.
\begin{proposition}
If $\{\bold{B}_n,\bold{C}_n\}$ is a minimizing sequence for $\mathfrak{J}_{red}$, then
$\{ \tilde{\bold{A}}[\bold{B}_n,\bold{C}_n],\bold{B}_n,\bold{C}_n\}$ is a 
minimizing sequence of $\mathfrak{J}$ and 
the equality, 
\[ \inf_{\bold{A},\bold{B},\bold{C}} \mathfrak{J}(\bold{A},\bold{B},\bold{C}) = 
\inf_{\bold{B},\bold{C}} \mathfrak{J}_{red}(\bold{B},\bold{C}), \]
holds.
\end{proposition}
\begin{proof}
Given that $\{\bold{B}_n,\bold{C}_n\}$ 
is a minimizing sequence of $\mathfrak{J}_{red}$:
$\lim_{n\to \infty} \mathfrak{J}_{red}(\bold{B}_n,\bold{C}_n) 
\to \inf_{\bold{B},\bold{C}} \mathfrak{J}_{red}(\bold{B},\bold{C})$. Since 
$\inf_{\bold{A},\bold{B},\bold{C}} \mathfrak{J}(\bold{A},\bold{B},\bold{C}) 
\leq \mathfrak{J}(\tilde{\bold{A}}[\bold{B}_n,\bold{C}_n],\bold{B}_n,\bold{C}_n) = 
\mathfrak{J}_{red}(\bold{B}_n,\bold{C}_n),$
we obtain
\[\inf_{\bold{A},\bold{B},\bold{C}} \mathfrak{J}(\bold{A},\bold{B},\bold{C}) 
 \leq \inf_{\bold{B},\bold{C}} \mathfrak{J}_{red}(\bold{B},\bold{C})\]
by passing to the limit. On the other hand, $ \mathfrak{J}(\bold{A},\bold{B},\bold{C}) 
\geq \mathfrak{J}(\tilde{\bold{A}}[\bold{B},\bold{C}],\bold{B},\bold{C}) = 
\mathfrak{J}_{red}(\bold{B},\bold{C}) \geq \inf_{\bold{B},\bold{C}} \mathfrak{J}_{red}(\bold{B},\bold{C})$
for arbitrary $\bold{A},\bold{B},\bold{C}$.
It follows that $\inf_{\bold{A},\bold{B},\bold{C}} \mathfrak{J}(\bold{A},\bold{B},\bold{C}) \geq 
\inf_{\bold{B},\bold{C}} \mathfrak{J}_{red}(\bold{B},\bold{C})$.
\end{proof}

\begin{cor}
If $(\bold{B}_*,\bold{C}_*)$ are minimizers of $\mathfrak{J}_{red}$, then 
$(\tilde{\bold{A}}[\bold{B}_*,\bold{C}_*],\bold{B}_*,\bold{C}_*)$ are minimizers of $\mathfrak{J}$. 
\end{cor}
\begin{proof}
From Lemma \ref{fixedBC}, the factor $\tilde{\bold{A}}[\bold{B},\bold{C}]$ always exists. Then 
if a minimizer $(\bold{B}_*,\bold{C}_*)$ of $\mathfrak{J}_{red}(\bold{B},\bold{C})$ exists, 
then $(\tilde{\bold{A}}(\bold{B}_*,\bold{C}_*),\bold{B}_*,\bold{C}_*)$ also
exists and it is a minimizer of $\mathfrak{J}$.
\end{proof}

\subsection{Analysis of the reduced objective function}
The introduction of $\mathfrak{J}_{red}$ reduces the number of unknown factors by one. 
In this section, we explicitly calculate $\mathfrak{J}_{red}$. Define
\begin{eqnarray}\label{bigM}
 \mathcal{M}_{\alpha\beta\gamma\delta}:= \sum_{i=1}
 \mathcal{T}_{i\alpha\beta} \mathcal{T}_{i\gamma\delta} \in \mathbb{R}^{J \times K \times J \times K},
 \end{eqnarray}
 a fourth order tensor from a contracted product over one index of two identical third-order tensors. 
The matricization $\bold{M} \in \mathbb{R}^{JK \times JK}$ of 
$\mathcal{M} \in \mathbb{R}^{J \times K \times J \times K}$ is defined by the following: 
\[ (\mathcal{M})_{\alpha\beta\gamma\delta} \longrightarrow (\bold{M})_{ij} \]
where $i=[\alpha + (\beta -1)J]$ and $j=[\gamma + (\delta-1)J]$. From (\ref{bigM}), 
we have the symmetry $(\mathcal{M})_{\alpha\beta\gamma\delta}=(\mathcal{M})_{\gamma\delta\alpha\beta}$ 
which implies that the matrix $\bold{M}$ is symmetric; i.e. $\bold{M}_{ij}=\bold{M}_{ji}$. 
It was shown in \cite{BLNT} that due to the isomorphic group 
structures between the sets of invertible tensors and matrices: 
a symmetric $((\mathcal{M})_{uvwx}=(\mathcal{M})_{wxuv})$ fourth order
 tensor $\mathcal{M}^{J \times K \times J \times K}$ has an eigendecomposition:
\begin{eqnarray}\label{fourthorderEVD}
\mathcal{M} = \bar{\mathcal{V}} \ast \mathcal{S} \ast \bar{\mathcal{V}}^T 
\end{eqnarray}
where $\ast$ is the contracted product of fourth order tensors defined as 
$(\mathcal{A} \ast \mathcal{B})_{ij\hat{i}\hat{j}} = \sum_{kl} (\mathcal{A})_{ijkl} (\mathcal{B})_{kl\hat{i}\hat{j}}$ given that the symmetric matrix $\bold{M}$ has an eigendecomposition such as $\bold{M}=\bar{\bold{V}} \cdot \bold{S} \cdot \bar{\bold{V}}^T$ where $(\bold{M})_{lm} \rightarrow \mathcal{M}_{uvwx}$, $(\bar{\bold{V}})_{lm} \rightarrow (\bar{\mathcal{V}})_{uvwx}$ and $(\bold{S})_{lm} \rightarrow (\mathcal{S})_{uvwx}$ with
$u,w=1,\ldots,J$, $v=\frac{l-u + J}{J}$ and  $x=\frac{m-w + J}{J}$. 
Note that $\bar{\bold{V}}$ is an orthogonal matrix and $\bold{S}$ is a diagonal matrix.

In accordance with the notation of Section~\ref{sec:prel} we can state some useful
tensor-vector and tensor-matrix operations for $\bold{a}, \bold{c} \in \R^{J}$ and
 $\bold{b}, \bold{d} \in \R^{K}$: 
\begin{enumerate}
\item
 $\mathcal{M} \bullet_{1,2,3,4}(\bold{a},\bold{b},\bold{c},\bold{d}) := 
\sum_{\alpha,\beta,\gamma,\delta} \mathcal{M}_{\alpha\beta\gamma\delta} 
a_{\alpha} b_{\beta}  c_{\gamma} d _{\delta} \in \R$ 
\item
$(\mathcal{M} \bullet_{2,3,4}(\bold{b},\bold{c},\bold{d}))_{\alpha} =  
\sum_{\beta,\gamma,\delta} \mathcal{M}_{\alpha \beta\gamma\delta} b_{\beta}  c_{\gamma} d _{\delta} 
\in \R^I $
\item 
$(\mathcal{M} \bullet_{2,4}(\bold{b},\bold{d}))_{\alpha,\gamma}
= \sum_{\beta,\gamma,\delta} \mathcal{M}_{\alpha \beta\gamma\delta} b_{\beta}  d _{\delta} 
\in \R^{I \times I}$
\end{enumerate}
Observe that 
\begin{eqnarray}\label{bigMmat}
\mathcal{M} \bullet_{1,2,3,4}(\bold{a},\bold{b},\bold{c},\bold{d}) 
= (\mathcal{T} \bullet_{2,3}(\bold{a},\bold{b}))^T
(\mathcal{T}\bullet_{2,3}(\bold{c},\bold{d}))= (\bold{a} \otimes \bold{b})^T \bold{M} 
(\bold{c} \otimes \bold{d}).
\end{eqnarray}
\begin{proposition}
The following is the reduced objective:
\begin{eqnarray}\label{jjj}
\mathfrak{J}_{red}(\bold{B},\bold{C}) = 
\frac{1}{2} \left (  \|\mathcal{T}\|_F^2 -  
\sum_{r,s}^R   (\bold{G}^\dagger)_{sr} \mathcal{M} 
\bullet_{1,2,3,4}(\bold{b_r},\bold{c_r},\bold{b_s},\bold{c_s}) \right )
\end{eqnarray}
where $\mathcal{M} \in \mathcal{R}^{J \times K \times J \times K}$ is defined in 
(\ref{bigM}) and  $\bold{G}^\dagger$ is the the pseudo-inverse of $\bold{G}$  in (\ref{theG}).
\end{proposition}
\begin{proof}
Expanding (\ref{mainobj2}) yields
\[ \mathfrak{J}_{red}(\bold{B},\bold{C}) =\mathfrak{J}(\tilde{\bold{A}}[\bold{B},\bold{C}],\bold{B},\bold{C}) =
 \frac{1}{2} \left ( \langle \mathcal{T}, \mathcal{T} \rangle -
2 \left \langle  \mathcal{T}, \sum_r \bold{\tilde{a}_r} \circ \bold{b_r} \circ 
\bold{c_r} \right \rangle +  \left \langle \sum_r \bold{\tilde{a}_r} 
\circ \bold{b_r} \circ \bold{c_r} , \sum_r \bold{\tilde{a}_r} 
\circ \bold{b_r} \circ \bold{c_r}  \right \rangle \right). \]
Using the component-wise definition of  $\tilde{\bold{A}}[\bold{B},\bold{C}]$ in Lemma \ref{fixedBC}, 
\[ {\tilde{a}}_{ir} = 
\sum_{s=1}^R (\mathcal{T} \bullet_{2,3}(\bold{B},\bold{C}))_{is}(\bold{G}^{\dagger})_{sr}
= \sum_{s=1}^R  \sum_{jk}(\mathcal{T}_{ijk}b_{js}c_{ks})(\bold{G}^{\dagger}_{sr}),\]
we obtain the following:\\
1. \begin{eqnarray}
\left \langle  \mathcal{T}, \bold{\tilde{a}_r} 
\circ \bold{b_r} \circ \bold{c_r} \right \rangle&=&
\sum_{r=1}^R \sum_{i,j,k} \mathcal{T}_{ijk} b_{jr} c_{kr} \tilde{a}_{ir} 
 = \sum_{r,s}^R \left ( \sum_{\substack{  i \\ j,k,\bar{j},\bar{k}  } }  
\mathcal{T}_{ijk} b_{jr} c_{kr} 
\mathcal{T}_{i\bar{j}\bar{k}} b_{\bar{j}s} c_{\bar{k}s} \right ) (\bold{G}^{\dagger})_{sr} \\ \nonumber
&=& \sum_{r,s}^R   
(\bold{G}^\dagger)_{s,r} \mathcal{M} \bullet_{1,2,3,4}(\bold{b_r},\bold{c_r},\bold{b_s},\bold{c_s}) \label{term2}
\end{eqnarray}
2. \begin{eqnarray}
\left \langle \sum_r \bold{\tilde{a}_r} 
\circ \bold{b_r} \circ \bold{c_r}, \sum_r \bold{\tilde{a}_r} 
\circ \bold{b_r} \circ \bold{c_r} \right \rangle& =& \sum_{r,\hat{r}} \sum_{ijk} (\tilde{a}_{ir}b_{jr}c_{kr}) (\tilde{a}_{i\hat{r}}b_{j\hat{r}}c_{k\hat{r}}) 
= \sum_{r,\hat{r}} \left ( \sum_{i} \tilde{a}_{ir} 
\tilde{a}_{i\hat{r}} \sum_{j} b_{jr} b_{j\hat{r}}  \sum_{k} c_{kr}c_{k\hat{r}} \right ) \nonumber\\
& =&  \sum_{r,\hat{r}}  \langle \bold{\tilde{a}_r},\bold{\tilde{a}_{\hat{r}}} \rangle 
\bold{G}_{r\hat{r}}  \label{starteq}  \\
&=&  \sum_{r,\hat{r}}  \left \langle \sum_{s=1}^R  
\sum_{ijk}(\mathcal{T}_{ijk}b_{js}c_{ks})(\bold{G}^{\dagger})_{sr} , 
\sum_{\hat{s}=1}^R  \sum_{\hat{i}\hat{j}\hat{k}}
(\mathcal{T}_{\hat{i}\hat{j}\hat{k}}b_{\hat{j}\hat{s}}c_{\hat{k}\hat{s}})(\bold{G}^{\dagger})_{\hat{s}\hat{r}}
 \right \rangle \bold{G}_{r\hat{r}} \nonumber \\
&=& \sum_{r,\hat{r}} \sum_{i}  \left ( \sum_{s,\hat{s}}  \sum_{jk,\hat{j}\hat{k}} 
 (\mathcal{T}_{ijk}b_{js}c_{ks})(\bold{G}^{\dagger})_{sr}
 (\mathcal{T}_{i\hat{j}\hat{k}}b_{\hat{j}\hat{s}}c_{\hat{k}\hat{s}})(\bold{G}^{\dagger})_{\hat{s}\hat{r}}
 \right ) \bold{G}_{r\hat{r}} \nonumber 
\end{eqnarray}
From the Moore-Penrose properties, we find 
\begin{eqnarray}\label{term3}
\left \langle \sum_r \bold{\tilde{a}_r} 
\circ \bold{b_r} \circ \bold{c_r}, \sum_r \bold{\tilde{a}_r} 
\circ \bold{b_r} \circ \bold{c_r} \right \rangle
& =& \sum_{s,\hat{s}} (\bold{G}^\dagger \bold{G} \bold{G}^\dagger)_{\hat{s}s} 
  \mathcal{M} \bullet_{1,2,3,4}(\bold{b_s},\bold{c_s},\bold{b_{\hat{s}}},\bold{c_{\hat{s}}}) \nonumber \\
& =& \sum_{s,\hat{s}} \bold{G}^\dagger_{\hat{s}s}  
 \mathcal{M} \bullet_{1,2,3,4}(\bold{b_s},\bold{c_s},\bold{b_{\hat{s}}},\bold{c_{\hat{s}}})
\end{eqnarray}
Equations (\ref{term2}) and (\ref{term3}) imply that the reduced objective is given by (\ref{jjj}).
\end{proof}
We can further simplify the functional:
\begin{lemma} \label{simplifiedlemma}
Let $\bold{U},\bold{\Sigma},\bold{V}$ be the matrices in the singular value decomposition
of $(\bold{B} \odot \bold{C})$, i.e., 
$(\bold{B} \odot \bold{C})=\bold{U} \cdot \bold{\Sigma} \cdot \bold{V}^T \in \mathbb{R}^{JK \times R}$ with 
$\bold{U} \in \mathbb{R}^{JK \times JK}$ orthogonal, 
$\bold{\Sigma} \in \mathbb{R}^{JK \times R}$ diagonal and $\bold{V} \in \mathbb{R}^{R \times R}$ orthogonal. 
Then,
\[ \mathfrak{J}_{red}(\bold{B},\bold{C}) = 
\frac{1}{2} \left ( \|\mathcal{T}\|^2_F - 
\sum_{r=1}^{\bar{R}} \langle \bold{u_k},\bold{M} \bold{u_k} \rangle \right ) \]  
where $\bold{M}$ is the matricization of $\mathcal{M}$ in (\ref{bigM}) and 
$\bar{R}=\mbox{rank}(\bold{\Sigma}) = \mbox{rank}(\bold{B} \odot \bold{C}) $ and
$\bold{u_k}$ is the $k$-th column of $\bold{U}$. 
\end{lemma}
\begin{proof}
Starting from (\ref{term3}), with the shortcut $\bold{\tilde{A}} = \bold{\tilde{A}}[\bold{B},\bold{C}]$
and symmetry of $\bold{G}^\dagger$ we find
\begin{eqnarray*}
\sum_{r,s}^R   (\bold{G}^\dagger)_{sr} 
\mathcal{M} \bullet_{1,2,3,4}(\bold{b_r},\bold{c_r},\bold{b_s},\bold{c_s})&=&
\langle (\bold{B} \odot \bold{C}) \cdot \tilde{\bold{A}}^T,   
(\bold{B} \odot \bold{C}) \cdot \tilde{\bold{A}}^T \rangle
= 
\langle (\bold{B} \odot \bold{C}),    
(\bold{B} \odot \bold{C}) \cdot \tilde{\bold{A}}^T\cdot  \tilde{\bold{A}}\rangle
\\
&=& \langle (\bold{B} \odot \bold{C}) , 
 (\bold{B} \odot \bold{C}) \cdot \bold{G}^{\dagger} \cdot 
 \mathcal{T} \bullet_{2,3}(\bold{B},\bold{C})^T \mathcal{T} \bullet_{2,3}(\bold{B},\bold{C}) 
\cdot \bold{G}^{\dagger}   \rangle \\
%
%
%
%
&=& \langle (\bold{B} \odot \bold{C}) \cdot \bold{G}^{\dagger}, 
 (\bold{B} \odot \bold{C}) \cdot \bold{G}^{\dagger} \cdot 
 \mathcal{T} \bullet_{2,3}(\bold{B},\bold{C})^T \mathcal{T} \bullet_{2,3}(\bold{B},\bold{C})   \rangle \\
&=& \langle  (\bold{B} \odot \bold{C}) \cdot \bold{G}^{\dagger},  
(\bold{B} \odot \bold{C}) \cdot \bold{G}^{\dagger} \cdot (\bold{B} 
\odot \bold{C})^T \bold{M} (\bold{B} \odot \bold{C})  \rangle~~~~~~(\mbox{from~}\ref{bigMmat})\\
&=& \langle  (\bold{B} \odot \bold{C})\cdot  \bold{G}^{\dagger} \cdot 
(\bold{B} \odot \bold{C})^T ,  (\bold{B} \odot \bold{C}) \cdot \bold{G}^{\dagger} \cdot 
(\bold{B} \odot \bold{C})^T \cdot \bold{M}   \rangle\\
&=& \mbox{Tr}(\bold{B} \odot \bold{C}) \cdot \bold{G}^{\dagger} \cdot 
 (\bold{B} \odot \bold{C})^T(\bold{B} \odot \bold{C}) \cdot \bold{G}^{\dagger} \cdot 
(\bold{B} \odot \bold{C})^T \cdot \bold{M} \\
 &=&\mbox{Tr}\bold{P}_{\bold{B}\odot \bold{C}} \cdot \bold{M}, 
\end{eqnarray*}
where $\mbox{Tr}$ denotes the matrix trace and 
\begin{eqnarray}\label{projector}
\bold{P}_{\bold{B}\odot \bold{C}}= 
(\bold{B} \odot \bold{C}) \cdot  [ (\bold{B} \odot \bold{C})^T \cdot  
(\bold{B} \odot \bold{C}) ]^{\dagger} \cdot
(\bold{B} \odot \bold{C})^T \cdot (\bold{B} \odot \bold{C}) \cdot
 [ (\bold{B} \odot \bold{C})^T \cdot (\bold{B} \odot \bold{C}) ]^{\dagger} \cdot (\bold{B} \odot \bold{C})^T.
\end{eqnarray}

Since $(\bold{B} \odot \bold{C})=\bold{U} \cdot \bold{\Sigma} \cdot  \bold{V}^T$ and
 $\bold{G}^{\dagger}=\bold{V} \cdot (\bold{\Sigma}^T \cdot \bold{\Sigma})^{\dagger} \cdot \bold{V}^T$, 
it holds that
\begin{eqnarray*}
\sum_{r,s}^R   (\bold{G}^\dagger)_{sr} \mathcal{M} \bullet_{1,2,3,4}(b_r,c_r,b_s,c_s)= 
\mbox{Tr}\bold{U} \cdot \bold{P}_{\Sigma} \cdot \bold{P}_{\Sigma}^T \cdot \bold{U}^T \cdot \bold{M}
=\mbox{Tr}(\bold{U}\cdot \bold{P}_{\Sigma})^T \cdot \bold{M}  \cdot \bold{U} \cdot \bold{P}_{\Sigma}
\end{eqnarray*}
where the projector matrix   $\bold{P}_{\Sigma}=\bold{\Sigma} (\bold{\Sigma}^T\bold{\Sigma})^{\dagger} 
\bold{\Sigma}^T \in \mathbb{R}^{JK \times JK}$ can be calculated as 
\[ (\bold{P}_{\Sigma})_{ij} = \left\{ \begin{array}{cc} 1 & \mbox{ if } i = j \mbox{ and }  i\leq \mbox{rank}(\bold{\Sigma})  \\
                          0 & \mbox{else} 
                         \end{array} \right.\]
Thus, finally 
\[ \sum_{r,s}^R   (\bold{G}^\dagger)_{sr} \mathcal{M} \bullet_{1,2,3,4}(b_r,c_r,b_s,c_s)=
\sum_{r=1}^{\bar{R}} \langle \bold{u_r},\bold{M} \bold{u_r} \rangle.\]                        
\end{proof}

The previous lemma allows us to rewrite the minimization problem 
for $ \mathfrak{J}_{red}$ into a Rayleigh quotient type problem.
\begin{thm}\label{thmyprop}
The minimization problem for $\mathfrak{J}_{red}(\bold{B},\bold{C})$
is equivalent to the following maximization problem
\begin{equation}\label{maxprob}
    \sup_{\bold{u_1},\ldots, \bold{u_{\bar{R}}}} \sum_{r=1}^{\bar{R}} \langle
 \bold{u_r},\bold{M} \bold{u_r} \rangle,
\end{equation}
where $\bold{u_1},\ldots, \bold{u_{\bar{R}}}$  is an orthonormal basis of 
$\mbox{range}(\bold{B} \odot \bold{C})$ with $\bar{R}=rank(\bold{B} \odot \bold{C})$. 
Equivalence holds in the following sense: if $(\bold{B},\bold{C})$ are (approximate) minimizers of
 $\mathfrak{J}_{red}(\bold{B},\bold{C})$, then any orthonormal basis of  
$\mbox{range}(\bold{B} \odot \bold{C})$
is a(n) (approximate) maximizer of \eqref{maxprob}. Conversely, if $\bold{u_1},\ldots \bold{u_{\bar{R}}}$
 are (approximate)
maximizers of \eqref{maxprob}, then the associated $(\bold{B},\bold{C})$ are (approximate) minimizers of
 $\mathfrak{J}_{red}(\bold{B},\bold{C})$.
\end{thm}
\begin{proof}
From Lemma~\ref{simplifiedlemma} it is clear that (approximate) minimizers of  $\mathfrak{J}_{red}$ 
are equivalent to (approximate) maximizers (\ref{maxprob}) over the left singular 
vectors of $\bold{B} \odot \bold{C}$. The maximization in (\ref{maxprob}) can be equally well done
over any orthonormal basis of the range of $\bold{B} \odot \bold{C}$: let $\bold{U}$ be the  
$\R^{IK \times \bar{R}}$
matrix with columns the left singular vectors corresponding to nonzero singular values. The column vectors
are an orthonormal basis of $\bold{B} \odot \bold{C}$. Similarly, for any other orthonormal basis 
of this range we can build a matrix $\bold{W}$ with columns the basis vectors, which is related
to $\bold{U}$ by $\bold{W} = \bold{U} \bold{Q}$, where $\bold{Q}$ is an $\bar{R}\times \bar{R}$ orthonormal matrix. 
By the invariance of the trace, the sum in (\ref{maxprob}) can be written as
\[ \sum_{r=1}^{\bar{R}} \langle
 \bold{u_r},\bold{M} \bold{u_r} \rangle,
= \mbox{Tr}(\bold{U}^T \bold{M} \bold{U}) = 
\mbox{Tr}(\bold{Q}^T\bold{U}^T \bold{M} \bold{U} \bold{Q}) =
 \mbox{Tr}(\bold{W}^T \bold{M} \bold{W}), \]
which ends the proof. 
\end{proof}

The new problem  formulation  (\ref{maxprob})  clearly indicates why the least squares problem might not
have a solution. Obviously, the functional $\sum_{r=1}^{\bar{R}} \langle \bold{u_r},\bold{M} \bold{u_r} \rangle,$
is continuous with respect to $\bold{u_1},\ldots \bold{u_{\bar{R}}}$ and since these vectors are orthonormalized
they are within a compact set. However, the additional restriction that
$\{\bold{u_1},\ldots \bold{u_{\bar{R}}}\}$ spans the range of a Khatri-Rao product space does not necessarily induce a closed 
set. In fact, a problem arises when the rank of the Khatri-Rao product decreases for 
a minimizing sequence.

\begin{proposition}\label{myprop}
Let $(\bold{B}_n,\bold{C}_n)$ be a minimizing sequence of  $\mathfrak{J}_{red}$ (and thus
$(\bold{A}[\bold{B}_n,\bold{C}_n],\bold{B}_n,\bold{C}_n)$ a minimizing sequence of
$\mathfrak{J}$). Without loss of generality, we can assume that there exists matrices
$\tilde{\bold{B}}$ and $\tilde{\bold{C}}$ with
\begin{equation}\label{eqhelp}
 \lim_{n\to \infty} \bold{B}_n = \tilde{\bold{B}} \quad \mbox{~~and~~}
 \lim_{n\to \infty} \bold{C}_n= \tilde{\bold{C}}.
\end{equation}
If the following rank condition,
\begin{equation}\label{rankcond} 
\liminf_{n\to \infty} \mbox{rank} (\bold{B}_n  \odot \bold{C}_n )
 \leq   \mbox{rank} (\bold{B} \odot \bold{C} ).
\end{equation}
holds, then $(\tilde{\bold{B}},\tilde{\bold{C}})$ is a minimizer of $\mathfrak{J}_{red}$ and
$(\bold{A}[\tilde{\bold{B}},\tilde{\bold{C}}],\tilde{\bold{B}},\tilde{\bold{C}})$ 
is a minimizer of $\mathfrak{J}$. In particular,
a solution to the least squares problem exists.  
\end{proposition}
\begin{proof}
We first notice that $\mathfrak{J}_{red}(\bold{B}_n,\bold{C}_n)$ does not change when 
 $(\bold{B}_n,\bold{C}_n)$ is replaced by 
$\frac{\bold{B}_n}{\|\bold{B}_n,\|},\frac{\bold{C}_n}{\|\bold{C}_n,\|}$ since the range 
of the Khatri-Rao product does not change. 
So if the original sequence is a minimizing sequence, then so is
 $(\frac{\bold{B}_n}{\|\bold{B}_n,\|},\frac{\bold{C}_n}{\|\bold{C}_n,\|})_n$.
By a compactness argument,  these matrices have a converging subsequence, which must again be a minimizing 
sequence. Thus, there is no loss of generality in assuming that the minimizing sequence of matrices
converges as in \eqref{eqhelp}. Let $r_n$ be  the rank of $(\bold{B}_n \odot \bold{C}_n)$ such that
$r^* = \liminf_{n} r_n$
and $r=\mbox{rank}(\bold{B} \odot \bold{C})$. By using a subsequence argument, we can
assume without loss of generality that $\lim_n {r_n} \to r^*$.
Now let us consider the associated left singular vectors of
 $(\bold{B}_n \odot \bold{C}_n)$. 
The sequence of vectors $(\bold{u_1^n},\ldots,\bold{u_{JK}^n)}$ are normalized eigenvectors of 
$(\bold{B}_n \odot \bold{C}_n) (\bold{B}_n \odot \bold{C}_n)^T$ and by compactness
 we can find another subsequence for which
all eigenvalues of $(\bold{B}_n \odot \bold{C}_n)(\bold{B}_n \odot \bold{C}_n)^T$ converge:
\[ \bold{u_i^n} \to_{n\to \infty} \bold{w_i} \quad i = 1,\ldots JK\]
If $\bold{u_i^n}$ corresponds to an eigenvalue $\lambda_i^n = 0$, for $n$ sufficiently large, it is obvious 
that $\bold{w_i}$ is in the nullspace of  $(\bold{B} \odot \bold{C}) (\bold{B} \odot \bold{C})^T$. On the other hand, if 
$u_i^n$ corresponds to an eigenvalue with  $\liminf_n \lambda_i^n >0$, then since the eigenvalues are
continuous functions of the matrix we get for a subsequence that
\[ \bold{w_i} = \lim_{n\to \infty} \bold{u_i^n} = 
\lim_{n\to \infty}  \frac{1}{\lambda_i^n}(\bold{B}_n \odot \bold{C}_n) 
(\bold{B}_n \odot \bold{C}_n)^T \bold{u_i^n} = 
 \frac{1}{\lambda_i} (\bold{B} \odot \bold{C})(\bold{B} \odot \bold{C})^T \bold{w_i},
\]
thus, $\bold{w_i}$ is also an eigenvalue of $(\bold{B} \odot \bold{C}) (\bold{B} \odot \bold{C})^T $.
With $r_n \to r_*$ we obtain  
\[ \bold{u_i^n} \to \bold{w_i} \quad i = 1,\ldots r^*. \]
Let us  denote by $\bold{w_i}$ the remaining eigenvectors spanning the range of $(\bold{B} \odot \bold{C})$
and let $N$ be the supremum in \eqref{maxprob}. Then 
\[ N  = \lim_{n\to  \infty} \sum_{i=1}^{r_n} \langle \bold{u_i^n}, \bold{M} \bold{u_i^n} \rangle
= \sum_{i=1}^{r_*}  \langle \bold{w_i}, \bold{M} \bold{w_i} \rangle \leq 
 \sum_{i=1}^{r} \langle \bold{w_i}, \bold{M} \bold{w_i} \rangle  \leq N \]
which shows that equality holds in this formula and thus, 
$(\bold{w_i})_{i=1}^r$ are maximizers of \eqref{maxprob} and the associated matrices
$(\bold{B},\bold{C})$ are minimizers of $\mathfrak{J}_{red}$. 
%
\end{proof}

 
Converse to these propositions is the following result that if a minimizer does not exist, then 
the rank of the Khatri-Rao product must change in the limit for any minimizing sequence. 
More precisely, the rank of $(\bold{B} \odot \bold{C})$ of the limit of a minimizing
sequence must be lower than the limit of the rank of $(\bold{B}_n \odot \bold{C}_n)$. 
If this is the case,  at least one singular value of $(\bold{B}_n \odot \bold{C}_n)$ tends
to $0$. A consequence of this is that the pseudo-inverse $\bold{G}^\dagger$ becomes unbounded, and 
thus, the norm of $\bold{A}_n[\bold{B}_n,\bold{C}_n]$ may become unbounded.  
This reflect the well-known fact of diverging summands in, see, e.g., \cite{SilLim}, 
which is referred to as the degenerate CP case.

%
%

\subsection{Rank-1 approximation}
It is worthwhile to study the special case of a least squares approximation \eqref{cost1} with $R = 1$.
In this case, it is well-known that a minimizer always exists. Moreover, the minimizers can
be calculated by a Rayleigh quotient type maximization. From the previous calculations, we obtain the following:

\begin{corollary}
Consider the least squares problem 
\begin{eqnarray} \label{cost1}
\min_{\bold{a},\bold{b},\bold{c}}  \frac{1}{2} \Vert \mathcal{T} -  \bold{a} \circ \bold{b} \circ \bold{c} \Vert _F^{2}.
\end{eqnarray} 
Minimizers to this problem always exist and the vectors $\bold{b},\bold{c}$ can be found 
as the solution of either one of the following equivalent  problems
\begin{eqnarray}\label{thisequation}
& & \max_{\bold{b},\bold{c}} \frac{(\bold{b} \otimes \bold{c})^T \bold{M} 
(\bold{b} \otimes \bold{c})}{\|\bold{b}\|^2 \|\bold{c}\|^2} =
  \max_{\bold{b},\bold{c}} 
\frac{\mathcal{M} \bullet_{1,2,3,4}(\bold{b},\bold{c},\bold{b},\bold{c})}{\|\bold{b}\|^2 \|\bold{c}\|^2} \\
& &  = \max_{ \substack{  \|\bold{b}\|=1 \\ \|\bold{c}\|=1} } 
(\bold{b} \otimes \bold{c})^T \bold{M} (\bold{b} \otimes \bold{c})
= \max_{ \substack{  \|\bold{b}\|=1 \\ \|\bold{c}\|=1} } 
\mathcal{M} \bullet_{1,2,3,4}(\bold{b},\bold{c},\bold{b},\bold{c})
\end{eqnarray}
\end{corollary}
\begin{proof}
In the case $R =1$, the Khatri-Rao product $\bold{b} \odot \bold{c}$ reduces to 
$\bold{b} \otimes \bold{c}$. For any $\bold{b}, \bold{c}$, 
$\frac{\bold{b}}{\|\bold{b}\|} \otimes \frac{\bold{c}}{\|\bold{c}\|}$ 
yields a (one-dimensional) orthonormal basis 
of the range of $\bold{b} \odot \bold{c}$. On the other hand, any normalized basis 
(which contains only one vector)
can be written as a Kronecker product with normalized vectors $\|\bold{b}\|=1,\|\bold{c}\|=1$. 
Proposition~\ref{myprop} yields the
equivalence of these problems. Without loss of generality we can take a minimizing sequence 
$(\bold{b_n},\bold{c_n})$
normalized to one. Since then neither $\bold{b_n}$ nor $\bold{c_n}$ can be zero vectors 
 $\mbox{rank}(\bold{b}_n \odot \bold{c}_n) =1$, and the rank of possible limit vectors 
is $\mbox{rank}(\bold{b} \odot \bold{c}) =1$. Thus, \eqref{rankcond} holds and a minimizer always exists. 
\end{proof}
The maximizers  in this corollary corresponds to the generalized singular values of 
$\bold{M}$ which was already proven in \cite{HOOI} by De Lathauwer, De Moor and Vanderwalle. 
But such characterization only holds in the case $R=1$.  Proposition~\ref{myprop} 
gives the generalization to $R >1$.

The optimality condition for the generalized Rayleigh quotient is well-known:
\begin{lemma}
A necessary condition for a maximizers $\bold{b},\bold{c}$ with $\|\bold{b}\| = 1, 
\|\bold{c}\| = 1$ in \eqref{thisequation}
is that there exists a number $\lambda$ such that
\begin{equation}\label{myeq}
 \mathcal{M} \bullet_{2,3,4}(\bold{c},\bold{b},\bold{c}) = \lambda \bold{b} \quad \mbox{~and~~}
 \mathcal{M} \bullet_{1,3,4}(\bold{b},\bold{b},\bold{c}) = \lambda \bold{c}. 
\end{equation}
\end{lemma}
Another way of stating the optimality condition is that $\bold{b}$ is an eigenvector of
the matrix $ \mathcal{M} \bullet_{2,4}(\bold{c},\bold{c})$ and $\bold{c}$ is an eigenvector 
(with the same eigenvalue) for the matrix
$ \mathcal{M} \bullet_{1,3}(\bold{b},\bold{b})$.
Unfortunately, straightforward linear algebra techniques cannot be applied to this nonlinear problem
 of calculating the maximizers  since the eigenvectors are coupled to each other.

We now look in more detail to the critical points of the functionals \eqref{thisequation}.
By a compactness argument, it is clear that the functional
\begin{eqnarray}  f(\bold{b},\bold{c}) : S^{I-1} \times S^{J-1} &\to& \R \nonumber \\
 \bold{b},\bold{c} & \mapsto &  \mathcal{M} \bullet_{1,2,3,4}(\bold{b},\bold{c},\bold{b},\bold{c}), \label{myfunc}
\end{eqnarray}
always has a maximum and a minimum where $S^n = \{\bold{x} \in \R^{n+1}\,|\, \|\bold{x}\| = 1 \}$ 
denotes the $n$-dimensional sphere. Moreover, for each maximum and minimum $\bold{b},\bold{c}$, 
the corresponding antipodal points $(\bold{b},-\bold{c}),$ $(-\bold{b},\bold{c}),$ $(-\bold{b},-\bold{c})$ 
are maxima and minima as well. 
Topologically, however, there must exist  critical points of this functional which are 
neither maxima nor minima, except in degenerate cases.
 In the two-dimensional
cases $I=J=2$, a related result was shown by De Lathauwer et al. \cite{LMV00}. We have a general result
in arbitrary dimensions. 
\begin{proposition}
Suppose that the maxima and minima of \eqref{myfunc} are nondegenerate in the sense
that the Hessian of $f$ at these points is non-singular (and thus, the Hessian is either 
negative or positive definite). Then there exist at least
4 additional critical points $\bold{b},\bold{c}$ that are neither maxima or minima of \eqref{myfunc}.
If all critical points are nondegenerate, then the number of critical points with index $\gamma$,
$\gamma = 0,\ldots (I-1)+(J-1)$ must 
satisfy the following conditions
\[ C_0 \geq 4, C_{(I-1)+(J-1)} \geq 4 \quad C_\gamma \mbox{ is divisible by } 4 \]
and
\[ R_\gamma - R_{\gamma_1} + \ldots \pm R_0  \leq C_\gamma - C_{\gamma_1} + \ldots \pm C_0 
 \quad \forall \gamma = 0,\ldots (I-1)+(J-1)
\]
where the $R_\gamma$ (the Betti-numbers) are the coefficients in the polynomial 
\[ (1+ x^{I-1})(1+x^{J-1}) = \sum_{\gamma=0}^{(I-1)+(J-1)} R_\gamma x^\gamma. \]
%
\end{proposition}
\begin{proof}
Suppose that besides the maximal and minimal points there are no additional 
critical point. Then, by the nondegeneracy condition, $f$ is a Morse function \cite{Mi63}.
However, by the Morse inequalities this is impossible. In fact, we know that there exist
at least 4 points of maxima and 4 points of minima. They correspond to critical point with index
$\gamma = (I-1)+(J-1)$ and $\gamma = 0$. Hence, denoting by $C_\gamma$ number of  critical points
with index $\gamma$, we have $C_{(I-1)+(J-1)} \geq 4$ and $C_0 \geq 4$. On the other hand, 
the Poincar\'e polynomial \cite{Ha00} of $S^{I-1} \times S^{J-1}$ is $(1+ x^{I-1})(1+x^{J-1})
= 1 + x^{J-1} + x^{I-1} + x^{(I-1)+ (J-1)}$, by the Morse inequalities 
$C_{I-1} \geq 1$ and $C_{J-1}\geq 1$, which imply the existence of 
critical points (neither being a maximum or minimum) of index $I-1$ and $J-1$.
Since for a critical point the corresponding antipodal points will be critical as
well, we have shown the existence of at least 4 critical points.

If we assume a-priori that all critical points are nondegenerate, $f$ will be 
a Morse function and the Morse inequalities as stated in the proposition
must be satisfied. By the same antipodal-point argument, the number $C_\gamma$ must always be divisible by $4$. 
\end{proof}
For the case $I=J=2$ we obtain 
$C_0 \geq 4$, $C_2 \geq 4$, $C_1 -C_0 \geq 2$, which imply that $C_1 \geq 6$ and $C_1 \geq 8$ due to
divisibility by $4$. Thus, even if we consider antipodal points as being equivalent, there must be at  
least two more critical points beside the extrema. 
In the case $I = J = 3$ the inequalities yield lower bounds $C_1 \geq 4, C_2 \geq 4, C_3 \geq 4$.
Up to antipodal points we have here at least three more critical points occurring in the 
case that all critical points are nondegenerate.   

The critical points of $f(\bold{b},\bold{c})$ correspond to critical points of the original
 least squares functional:
\begin{lemma}
 Let $\bold{b},\bold{c} \in  S^{I-1} \times S^{J-1}$ be a critical point of \eqref{myfunc},
then with the setting (cf. \eqref{theA}) $\bold{a} = \mathcal{T} \bullet_{2,3} (\bold{b},\bold{c})$
the vectors $(\bold{a},\bold{b},\bold{c})$ satisfy the first order optimality conditions  of \eqref{cost1}.
\end{lemma}
\begin{proof}
With the definition of $\mathcal{M}$ and \eqref{myeq},
a critical point $(\bold{b},\bold{c})$ satisfies
\[ \lambda \bold{b} = \sum_{i} \mathcal{T}_{i,2,3}(\bold{b},\bold{c}) 
 \mathcal{T}_{i,\bullet,3}(\bold{c}) = 
\mathcal{T}_{1,\bullet,3}(\bold{a},\bold{c}) 
\]
and
\[  \lambda \bold{c} = \mathcal{T}_{1,2,\bullet}(\bold{a},\bold{b}) 
\]
 $\lambda$ in the optimality condition can be calculated to 
\[ \lambda = \mathcal{M}_{1,2,3,4}(\bold{b},\bold{c},\bold{b},\bold{c}) = 
 \sum_{i} \mathcal{T}_{i,2,3}(\bold{b},\bold{c})  \sum_{i} \mathcal{T}_{i,2,3}(\bold{b},\bold{c}) = 
\|\bold{a}\|^2
\]
Thus we obtain the optimality conditions for \eqref{cost1}:
\begin{eqnarray}
 \mu \bold{\tilde{a}} &=& \mathcal{T} \bullet_{2,3} (\bold{b},\bold{c}) \label{cp1} \\
 \mu \bold{b} &=&  \mathcal{T}_{1,\bullet,3}(\bold{\tilde{a}},\bold{c}) \label{cp2} \\
 \mu \bold{c} &=&  \mathcal{T}_{1,2,\bullet}(\bold{\tilde{a}},\bold{b}), \label{cp3}
\end{eqnarray}
with $\bold{\tilde{a}} = \frac{\bold{a}}{\|\bold{a}\|}$ and $\mu = \|\bold{a}\|$. 
\end{proof}
Since ALS works with the first order optimality condition, it will saturate 
at a critical point. Thus, we have the following
negative result:
\begin{theorem}
If the extrema of  \eqref{myfunc} are nondegenerate, then there always exists
a set of vectors $(\bold{a},\bold{b},\bold{c})$ which is neither 
a maximum nor a minimum of \eqref{cost1} for which the ALS sequence for \eqref{cost1} remains constant
 at this point
\begin{equation}\label{saturate} (\bold{a_{k+1}},\bold{b_{k+1}},\bold{c_{k+1}}) = 
(\bold{a_{k}},\bold{b_{k}},\bold{c_{k}}) \quad \forall k \geq 1.
\end{equation}
\end{theorem}
\begin{proof} 
Taking as starting point for the ALS iteration a critical point satisfying
\eqref{cp1}--\eqref{cp3}, with vectors $\bold{a_0},\bold{b_0},\bold{c_0}$ normalized to 
norm 1. The ALS iteration in the rank-1 case reads
\[ \bold{a_{k+1}} = \frac{ \mathcal{T}_{\bullet,2,3} 
(\bold{b_k},\bold{c_k})}{\|\bold{b_k}\|^2 \|\bold{c_k}\|^2} \quad 
  \bold{b_{k+1}} = \frac{\mathcal{T}_{1,\bullet,3}
 ( \bold{a_{k+1}},\bold{c_k})}{\| \bold{a_{k+1}}\|^2 \|\bold{c_k}\|^2}  \quad 
\bold{c_{k+1}} = \frac{ \mathcal{T}_{1,2,\bullet}
 ( \bold{a_{k+1}},\bold{b_{k+1}})}{\| \bold{a_{k+1}}\|^2 \| \bold{b_{k+1}}\|^2}
\]
It follows by induction that with the given starting value, the iteration becomes
\[ \bold{a_{k+1}} = \alpha_{k+1} \bold{a_0} \quad
 \bold{b_{k+1}} = \beta_{k+1} \bold{b_0} \quad 
\bold{c_{k+1}} = \gamma_{k+1} \bold{c_0},
\]
where $\alpha_{k+1},\beta_{k+1},\gamma_{k+1}$ are numbers satisfying the recursion
\[ \alpha_{k+1} = \frac{\mu}{\beta_k \gamma_k} \quad
 \beta_{k+1} = \frac{\mu}{\alpha_{k+1} \gamma_k} \quad
\gamma_{k+1} = \frac{\mu}{\alpha_{k+1} \beta_{k+1}} 
\]
for $k\geq 1$. Eliminating first $\alpha_{k+1}$ yields $\beta_{k+1} = \beta_{k}$, and furthermore
$\gamma_{k+1} = \gamma_k$ for all $k\geq 1$, hence $\alpha_{k+2} = \alpha_{k+1}$. 
Thus, we observe that the iteration remains constant \eqref{saturate}.  Since the extrema of $\mathfrak{J}$ are
one-to-one related to extrema of $\mathfrak{J}_{red}$ and hence of $f(\bold{b},\bold{c})$ 
the ALS sequence remains at a point which is not an extrema of the least squares functional.
\end{proof}
This result shows that there is no guarantee that a converging ALS sequence
yields a minimizer of $\mathfrak{J}$. Of course, this is not surprise for 
a first order method. 

\subsection{Reduced functional in projection form}
We now derive an alternative form of the reduced functional 
$\mathfrak{J}_{red}$ as a weighted distance to the Khatri-Rao space. 
This form will be useful in the next section to design a simple algorithm
for finding an initial guess to the minimization form. 

Based on Lemma~\ref{simplifiedlemma} we can simplify the reduced functional
taking into account the diagonalization of $\bold{M}$:
\begin{lemma}\label{lemmams}
Let $(\bold{B} \odot \bold{C})=\bold{U} \bold{\Sigma} \bold{V}^T \in
 \mathbb{R}^{JK \times R}$ and $\bold{M}=\bar{\bold{V}}\bold{S} \bar{\bold{V}}^T $ 
with orthogonal matrices: $\bold{U} \in \mathbb{R}^{JK \times JK}$, 
$\bold{V} \in \mathbb{R}^{R \times R}$, $\bar{\bold{V}} \in \mathbb{R}^{JK \times JK}$ 
and diagonal matrices: $\bold{\Sigma} \in \mathbb{R}^{JK \times R}$ and 
$\bold{S}\in \mathbb{R}^{JK \times JK}$ with $\bold{S} = \mbox{diag}(\lambda_i)$. 
Denote by $\bold{ \bar{v}_i}$ the columns of $\bar{\bold{V}}$, 
then,
\begin{equation}\label{aaaeq}  \mathfrak{J}_{red}(\bold{B},\bold{C}) =\frac{1}{2} \sum_{i=1}^{JK}
 \lambda_i \left(\|\bold{\bar{v}_i}\|^2 - \sum_{r=1}^{\bar{R}} \langle \bold{\bar{v}_i},
\bold{u_r} \rangle^2 \right)
=\frac{1}{2}\sum_{i=1}^{JK} \lambda_i \left(1 - \sum_{r=1}^{\bar{R}} \langle \bold{\bar{v}_i},
\bold{u_r} \rangle^2 \right)
\end{equation}
\end{lemma}
\begin{proof}
Observe that $\Vert \mathcal{T} \Vert^2_F=\sum_{jk}(\sum_i \mathcal{T}_{ijk}\mathcal{T}_{ijk})=
\sum_{jk} \mathcal{M}_{jkjk}=\mbox{trace}(\bold{M})=\mbox{Tr}(\bold{S}\bold{V}^T\bold{V})=
\sum_{i}^{JK} \lambda_i \Vert \bold{\bar{v}_i} \Vert^2_F=\sum_i \lambda_i$ 
since $\Vert \bold{\bar{v}_i} \Vert^2_F=1$.
With $\langle \bold{u_r},\bold{M} \bold{u_r} \rangle = 
\langle \bold{u_r}, \sum_{i}^{JK} \lambda_i \bold{\bar{v}_i} \bold{\bar{v}_i}^T   \bold{u_r} \rangle= 
\sum_i^{jk} \lambda_i \langle \bold{\bar{v}_i} ,\bold{u_r}\rangle^2$, the result follows.
\end{proof}

New we define the {\bf Khatri-Rao range}, i.e. the range of the matrix $\bold{B} \odot \bold{C}$.
This range is a subset of $\R^{IJ}$; for later use it is convenient to define the Khatri-Rao range 
by matricizing this range. As usual we denote the columns of the matrices $\bold{B}$ and $\bold{C}$ by
$\bold{b_i}$ and $\bold{c_i}$:   
\begin{equation}\label{krrange} \lin(\bold{B},\bold{C}):= 
\left \{\bold{X}=\sum_{i=1}^R \mu_i \bold{b_i} \circ \bold {c_i}  \in \mathbb{R}^{J \times K} \,|\,
 ~~where~~ \mu_i \in \R,  
 \right \} 
\end{equation}
It is obvious that $\bold{X} \in \R^{I\times J}$ is in the Khatri-Rao range  
$\bold{X} \in \lin(\bold{B},\bold{C})$
if and only if its vectorized version 
$\bold{X^{vec}} \in \R^{IK}$ is in the range of $\bold{B} \odot \bold{C}$
\[ \bold{X^{vec}} = (\bold{B} \odot \bold{C}) \mbox{ \boldmath{$\mu$}}. \]

%

\begin{theorem}\label{main}
Let $\bold{\bar{V}}_i \in \mathbb{R}^{J \times K}$ be the matricized version of
a vector $\bold{\bar{v}_i} \in \mathbb{R}^{JK}$ appearing in Lemma~\ref{lemmams}. 
With the notation of Lemma~\ref{lemmams}, the reduced least squares functional is simplified as
\begin{equation} \mathfrak{J}(\bold{B},\bold{C}) = 
\frac{1}{2} \sum_{i=1}^{JK} \lambda_i \left( \|\bar{\bold{V}}_i - \lin(\bold{B},\bold{C})\|_F \right)^2,
\end{equation}  
where $\|\bar{\bold{V}}_i - \lin(\bold{B},\bold{C})\|$ denotes the distance
 of $\bar{\bold{V}}_i$ to the linear subspace $\lin(\bold{B},\bold{C})$
\[ \|\bar{\bold{V}}_i - \lin(\bold{B},\bold{C})\|_F = 
\inf_{\bold{X} \in \lin(\bold{B},\bold{C})} \| \bar{\bold{V}}_i - \bold{X}\|_F. \]
\end{theorem}
\begin{proof}
Let $\bold{P}_{\bold{U}}$ be the orthogonal projector onto 
$\mbox{range}(\bold{B} \odot \bold{C})$ since the vectors $\bold{u_i}$ are an orthogonal basis of
this range we have 
$\bold{P}_{\bold{U}} \bold{\bar{v}_i}= \sum_{r=1}^{\bar{R}} \langle \bold{\bar{v}_i}, \bold{u_r} \rangle 
\bold{u_r}$.
The minimum distance between $\bold{\bar{v}_i}$ and $\mbox{range}(\bold{B} \odot \bold{C})$
 can be expressed by the projector as :
\begin{eqnarray*}
\inf_{\bold{\tilde{x}} \in \mbox{range}(\bold{B} \odot \bold{C})}
 \Vert  \bold{\bar{v}_i} - \bold{\tilde{x}}  \Vert^2_F = 
 \Vert  \bold{\bar{v}_i} - \bold{P}_{\bold{U}}\bold{\bar{v}_i}  \Vert^2_F.
\end{eqnarray*}
Moreover, $\Vert \bold{\bar{v}_i}   \Vert^2_F - 
\Vert  \bold{P}_{\bold{U}} \bold{\bar{v}_i}  \Vert^2_F= 
\Vert \bold{\bar{v}_i}   \Vert^2_F - \sum_r^{\bar{R}} \langle \bold{\bar{v}_i}, 
\bold{u_r} \rangle^2 = \Vert  \bold{\bar{v}_i} - \bold{P}_{\bold{U}}\bold{\bar{v}_i}  \Vert^2_F.$ 
Since inner products and norms are the same for a vector and its matricization, we 
obtain the result from \eqref{aaaeq}. 
%
\end{proof}

Observe that for a particular
 index $i$ there exists a set of indices $(\hat{j},\hat{k})$ such $i=\hat{j}- (\hat{k}-1)J$
 which implies that $(\bar{\bold{V}}_i)_{jk}$ is matrix representing the
 subtensor $\mathcal{\bar{V}}_{\hat{j}\hat{k}jk}$ (\ref{fourthorderEVD}).

A simple consequence of the previous theorem is the following. 
\begin{corollary}\label{swamp}
If $\bar{\bold{B}},\bar{\bold{C}}$ and $\bold{B},\bold{C}$ are matrices that span the same Khatri-Rao range  i.e.
\[ \lin( \bar{\bold{B}},\bar{\bold{C}}) = \lin(\bold{B},\bold{C}) \]
then
\[ \mathfrak{J}_{red}(\bar{\bold{B}},\bar{\bold{C}}) = \mathfrak{J}_{red}(\bold{B},\bold{C}). \]
\end{corollary}

\begin{remark}\label{remarkswamp}
If this corollary is applied to the case when $\mathfrak{J} = 0$ 
we obtain -- as a special case -- a uniqueness condition.  
The CP decomposition $(\bold{A},\bold{B},\bold{C})$ is called unique up to permutation and 
scaling if any alternative decomposition $(\bold{\bar{A}},\bold{\bar{B}},\bold{\bar{C}})$ 
satisfies $\bold{\bar{A}}=\bold{A} \bold{\Pi}\bold{\Lambda}_1$, 
$\bold{\bar{B}}=\bold{B} \bold{\Pi}\bold{\Lambda}_2$ and 
$\bold{\bar{C}}=\bold{C} \bold{\Pi}\bold{\Lambda}_3$ where
 $\bold{\Pi}$ is an $R \times R$ permutation matrix and $\bold{\Lambda}_j$ are 
nonsingular matrices such that $\prod_{j=1}^n \bold{\Lambda}_j=\bold{I}_R.$
%
%
Certainly,  if $\bold{\bar{B}}=\bold{B} \bold{\Pi}\bold{\Lambda}_2$ and 
$\bold{\bar{C}}=\bold{C} \bold{\Pi}\bold{\Lambda}_3$, 
then $\lin( \bar{\bold{B}},\bar{\bold{C}}) = \lin(\bold{B},\bold{C})$ and thus,
$\mathfrak{J}_{red}(\bar{\bold{B}},\bar{\bold{C}}) = \mathfrak{J}_{red}(\bold{B},\bold{C})$. 
From Corollary~\ref{swamp} we find that if a CP decomposition is unique up 
to scaling and permutation then 
$\lin( \bar{\bold{B}},\bar{\bold{C}}) = \lin(\bold{B},\bold{C})$ can only hold
when  $\bold{\bar{B}}$ and $\bold{\bar{C}}$ is a scaled and permuted version 
of $\bold{B}$ and $\bold{C}$. 
\end{remark}

\begin{remark}
The reduced functional and its analysis is equally well doable for higher order tensors 
as well, e.g., in a forth order decomposition
\[  (\mathcal{A})_{ijkl}= \sum_{r=1}^R a_{ir} b_{jr} c_{kr} d_{lr}. \]
The Khatri-Rao range $\lin(\bold{B},\bold{C})$ has to be replaced by the 
analogous set 
\[ \lin(\bold{B},\bold{C},\bold{D}) = \left \{
\mathcal{X}=\sum_{i=1}^R \mu_i \bold{b_i} \circ \bold {c_i} \circ \bold{d_i} 
\in \mathbb{R}^{J \times K \times L} \,|\,
 ~~where~~ \mu_i \in \R,  
 \right \} 
\]
\end{remark}

In our view, Corollary~\ref{swamp} displays  one possible reason for the swamping effect. We explain this
in the following subsection. 
\subsubsection{One explanation of swamping}
The swamping phenomenon describes the effect that iterations method for  
 minimizing the functional $\mathfrak{J}(\bold{A},\bold{B},\bold{C})$ 
exhibit  a long interval of iterations where the functional 
value remains almost constant and does not decrease. 
This is commonly seen in the ALS implementation.

From the definition of $\mathfrak{J}_{red}$ in (\ref{defJred}), 
\[ \mathfrak{J}(\bold{A}^{k+1},\bold{B}^k,\bold{C}^k) = \mathfrak{J}_{red}(\bold{B}^k,\bold{C}^k).  \]
Moreover, for an iteration of the alternating minimization (ALS) procedure, we obtain
\begin{eqnarray*}
\mathfrak{J}(\bold{A}^k,\bold{B}^k,\bold{C}^k) \geq \mathfrak{J}(\bold{A}^{k+1},\bold{B}^k,\bold{C}^k)
 =\mathfrak{J}_{red}(\bold{B}^k,\bold{C}^k) &\geq& \mathfrak{J}(\bold{A}^{k+1},\bold{B}^{k+1},\bold{C}^k) \\
&\geq& \mathfrak{J}(\bold{A}^{k+1},\bold{B}^{k+1},\bold{C}^{k+1})\\
& \geq& \mathfrak{J}_{red}(\bold{B}^{k+1},\bold{C}^{k+1}).
\end{eqnarray*}
Thus, the functional values  $\mathfrak{J}(\bold{A}^k,\bold{B}^k,\bold{C}^k)$ will
 behave in a similar way as $\mathfrak{J}_{red}(\bold{B}^k,\bold{C}^k)$. 

Corollary~\ref{swamp} can serve as one possible explanation of  the swamping effect.  
It shows, that the dependence of the least squares functional on the matrices $\bold{B},\bold{C}$ 
is rather low, as it only depends on the Khatri-Rao range $\lin(\bold{B},\bold{C})$. 
In particular, if $\lin(\bold{B}^k,\bold{C}^k) = \lin(\bold{B}^{k+j},\bold{C}^{k+j})$ 
for some $j$ iterations, then $\mathfrak{J}_{red}(\bold{B}^{k},\bold{C}^{k}) = 
\mathfrak{J}_{red}(\bold{B}^{k+j-1},\bold{C}^{k+j-1})$ and as a consequence
$\mathfrak{J}(\bold{A}^{k+j},\bold{B}^{k+j},\bold{C}^{k+j})$ will stay at the same value for these 
iterations.  
Moreover,  the set of matrices that span the same linear space can be quite 
large which explains the large region at which least squares functional attains the same value.
This also explains the increasing length of the swamps present in high-order $n\geq 4$ tensors; 
e.g., the subspace $\lin(\bold{B}^k,\bold{C}^k,\bold{D}^k)$ corresponding to the 
functional $ \mathfrak{J}(\bold{A}^{k},\bold{B}^{k},\bold{C}^{k},\bold{D}^k)=
\mathfrak{J}_{red}(\bold{B}^{k},\bold{C}^{k},\bold{D}^k)$ of a fourth-order 
tensor is spanned by a huge set of matrices of $\bold{B}$, $\bold{C}$ and $\bold{D}$. 
This reasoning can be  underpinned by numerical calculations.

\begin{figure}[htp]
  \begin{center}
   \subfigure[The plot on top depicts an ALS swamp while the bottom plot tracks the condition number of the matrix
$(\bold{B}^k \odot \bold{C}^k$~~$\bold{B}^{k+1} \odot \bold{C}^{k+1})$]{\label{fig:ALSCALSCP}
\includegraphics[width = 46 mm]{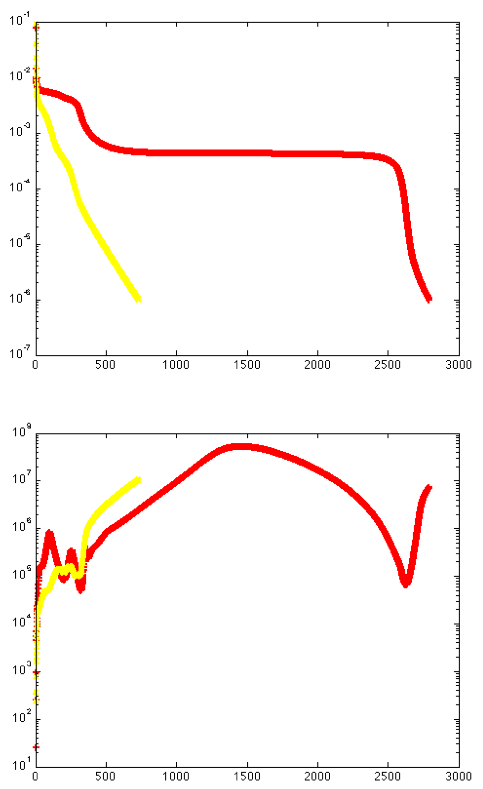}}
  \subfigure[Measurement of Subspaces. Top Left:
 $(\bold{B}^k \odot \bold{C}^k)$ vs $(\bold{B}^{k+1} \odot \bold{C}^{k+1})$, Middle Left:
 $(\bold{A}^k \odot \bold{C}^k)$ vs $(\bold{A}^{k+1} \odot \bold{C}^{k+1})$, Bottom Left:
 $(\bold{B}^k \odot \bold{A}^k)$ vs $(\bold{B}^{k+1} \odot \bold{A}^{k+1})$. 
Top Right: $\bold{A}^k$ vs $\bold{A}_{orig}$, Middle Right: $\bold{B}^k$ vs $\bold{B}_{orig}$,
 Bottom Right: $\bold{C}^k$ vs $\bold{C}_{orig}$.]{\label{fig:Subspaces}
\includegraphics[width = 71 mm]{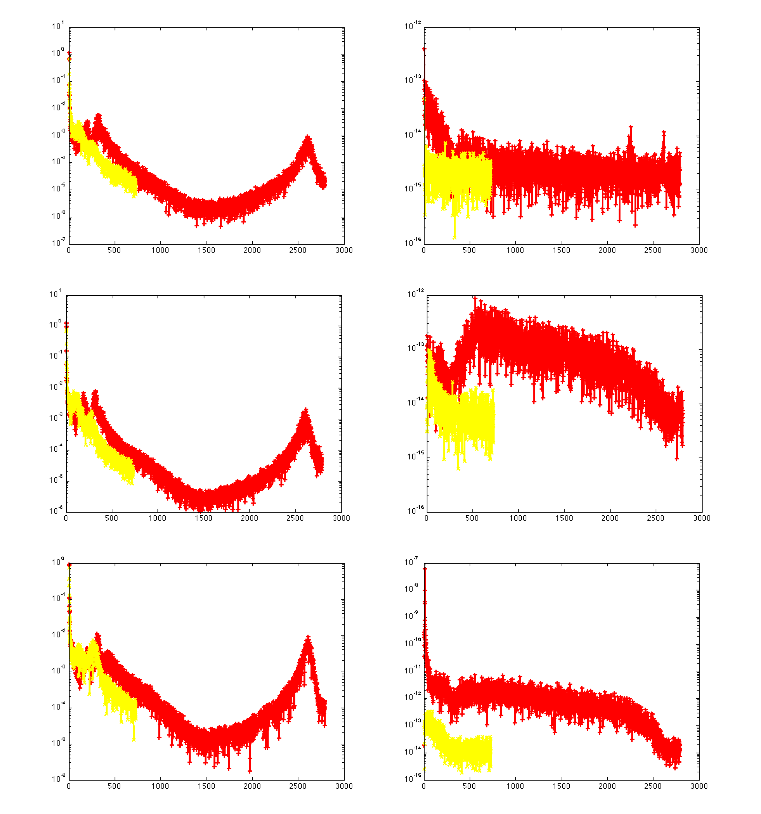}}
   \end{center}
  \caption{ALS (red -+-) and CALS (yellow -x-)}
  \label{fig:SwampySubspaces}
\end{figure}

In Figure~\ref{fig:SwampySubspaces}b, we measure the distance 
between subspaces $(\bold{B}^k \odot \bold{C}^k)$ and $(\bold{B}^{k+1} \odot \bold{C}^{k+1})$
 (top-left) by taking an arbitrary vector $\bold{x}$ and calculating the norm difference
 of the projections of $\bold{x}$ onto the spaces $(\bold{B}^k \odot \bold{C}^k)$ and 
$(\bold{B}^{k+1} \odot \bold{C}^{k+1})$.  As seen in Figure~\ref{fig:SwampySubspaces}b, 
at the swamp regime, the norm differences in the subspaces dip down to $10^{-6}$ in the 
ALS implementation which coincides with our swamp explaination that
$\lin(\bold{B}^k,\bold{C}^k) \approx \lin(\bold{B}^{k+j},\bold{C}^{k+j})$.  
The plots in Figure~\ref{fig:SwampySubspaces}b on the right column describe 
the measure of the subspaces spanned by $k$-th approximation $\bold{A}^k$ and the 
original factor $\bold{A}_{orig}$ (top-right). The norm distances are all 
fairly small, but relative to the yellow curves produced by using CALS (see Section~\ref{sec:numone}), $\bold{B}^k$ 
and $\bold{C}^k$ subspaces from ALS are far off from the original subspaces $\bold{B}_{orig}$
 and $\bold{C}_{orig}$ which, once again, indicates 
that $\lin(\bold{B}^k,\bold{C}^k) \approx \lin(\bold{B}^{k+j},\bold{C}^{k+j})$ for some $j$ 
accounting for the ALS swamp.

Another way to measure the distance between subspaces is through 
the condition number of the matrix $[(\bold{B}^k \odot \bold{C}^k)~~ (\bold{B}^{k+1} \odot \bold{C}^{k+1}) ]$
 as a way to measure linear independence. The bottom plot in  Figure~ \ref{fig:SwampySubspaces}a 
shows that when ALS is used, large condition numbers are present at the swamp regime, 
shooting up to $10^{9}$.

\section{Bounds on $\mathfrak{J}_{red}$ and suboptimal solutions}
In this section we prove lower and upper bound on $\mathfrak{J}_{red}$ using 
Theorem~\ref{main}. Moreover, we will also define a dominating functional 
$\mathfrak{L}$, which minimizers can be calculated by standard linear algebra methods.
The corresponding algorithm, called the Centroid Projection yields an initial guess, 
for a minimization step for $\mathfrak{J}$. 
 
From Theorem~\ref{main} we can use the Eckard-Young theorem to obtain lower bounds:
We keep the notation of Theorem~\ref{main} and Lemma~\ref{lemmams}:
\begin{corollary}\label{clb}
For all matrices $(\bold{B},\bold{C})$, the lower bound of $\mathfrak{J}_{red}$ is calculated as 
\begin{equation}
\inf_{(\bold{B},\bold{C})}  \mathfrak{J}_{red}(\bold{B},\bold{C}) \geq 
\frac{1}{2} \sum_{i=1}^{JK} \lambda_i  \left (\sum_{k=R+1}^{\min(J,K)} (\sigma_k^i)^2 \right )  
\end{equation}
if  $\sigma_{k}^i$ the $k$-th singular values of $\bold{\bar{V}}_i$. 
\end{corollary}
\begin{proof}
Recall that the Eckart-Young Theorem gives the infimum through the truncated SVD; i.e.
\[ \inf_{rank(\bold{X}) = R} \| \bar{\bold{V}}_i - \bold{X} \|_F=  
\left ( \sum_{k=R+1}^{\min\{J,K\}} (\sigma_{k}^i)^2 \right )^{\frac{1}{2}}\]
where $\sigma_{k}^i$'s are the singular value of $\bar{\bold{V}}_i$. 
Also, observe that $\lin(\bold{B},\bold{C})$ contains matrices with rank at most $R$, hence,
\begin{eqnarray*}
 \inf_{(\bold{B},\bold{C})}  \mathfrak{J}_{red}(\bold{B},\bold{C}) &=& 
\inf_{(\bold{B},\bold{C})}  \frac{1}{2} \sum_{i=1}^{JK} \lambda_i 
\inf_{\bold{X} \in \lin(\bold{B},\bold{C})} \left( \|\bar{\bold{V}}_i - \bold{X} \|_F \right)^2\\
 &\geq&   \frac{1}{2} \sum_{i=1}^{JK} \lambda_i 
\inf_{ \substack{ (\bold{B},\bold{C}) \\ \bold{X} \in (\bold{B},\bold{C}) }}
 \left( \|\bar{\bold{V}}_i - \bold{X} \|_F \right)^2\\
  &\geq&   \frac{1}{2} \sum_{i=1}^{JK} \lambda_i \inf_{ \mbox{rank}(\bold{X})
 \leq R} \left( \|\bar{\bold{V}}_i - \bold{X} \|_F \right)^2\\
&\geq&  \frac{1}{2} \sum_{i=1}^{JK} \lambda_i   \left (\sum_{k=R+1}^{\min(J,K)} (\sigma_k^i)^2 \right ). 
 \end{eqnarray*}
\end{proof}
This corollary can be used to find lower bounds on the distance 
of a tensor to its best rank $R$ approximation. In particular, if a tensor has rank $R$ 
it must hold that
\[ \sum_{i=1}^{JK} \lambda_i  \left (\sum_{k=R+1}^{\min(J,K)} (\sigma_k^i)^2 \right )   = 0.\] 
Note that the lower bound can a-priori be calculated by 
standard linear algebra method (eigenvalue and SVD decomposition).  
The computation requires an eigenvalue 
decomposition of $\bold{M}$ followed by a SVD of each of the matricized
 eigenvalues $\bold{\bar{V}}_i$. 

The next result establishes an upper bound using a dominating functional.
\begin{corollary}
For all matrices $(\bold{B},\bold{C})$, the upper bound of $\mathfrak{J}_{red}$ is calculated as 
\begin{equation}
\inf_{\bold{B},\bold{C}} \mathfrak{J}_{red}(\bold{B},\bold{C}) \leq
 \inf_{\bold{B},\bold{C}} \mathfrak{L}(\bold{B},\bold{C})  
\end{equation}
with
\begin{equation} \mathfrak{L}(\bold{B},\bold{C}):= 
 \inf_{\bold{X} \in \lin(\bold{B},\bold{C})} \sum_{i=1}^{JK} \lambda_i 
\left( \| \bold{\bar{V}}_i - \bold{X}\|_F^2 \right) 
= 
\end{equation}
Moreover,
\[ \inf_{\bold{B},\bold{C}} \mathfrak{L}(\bold{B},\bold{C})  = 
\inf_{ \mbox{\rm rank}(\bold{X}) \leq \bar{R}} \sum_{i=1}^{JK}
 \lambda_i \left( \|\bold{\bar{V}}_i  - \bold{X}\|_F^2 \right) \]
\end{corollary}
\begin{proof}
\begin{eqnarray*}
\inf_{(\bold{B},\bold{C})}  \mathfrak{J}_{red}(\bold{B},\bold{C}) &=& 
\inf_{(\bold{B},\bold{C})}  \frac{1}{2} \sum_{i=1}^{JK} \lambda_i \inf_{\bold{X} \in 
\lin(\bold{B},\bold{C})} \left( \|\bar{\bold{V}}_i - \bold{X} \|_F \right)^2\\
&\leq& \inf_{(\bold{B},\bold{C})} \inf_{\bold{X} \in \lin(\bold{B},\bold{C})}
 \frac{1}{2} \sum_{i=1}^{JK} \lambda_i \left( \|\bar{\bold{V}}_i - \bold{X} \|_F \right)^2\\
&=& \inf_{\mbox{rank}(\bold{X}) \leq R}  
\frac{1}{2} \sum_{i=1}^{JK} \lambda_i \left( \|\bar{\bold{V}}_i - \bold{X} \|_F \right)^2
\end{eqnarray*}
The last equality follows from the fact that $\lin(\bold{B},\bold{C})$
contains matrices with $\mbox{rank} \leq R$. Moreover using the SVD, for any matrix $\bold{S}$ of
rank at most $R$ it follows that  $\bold{S} \in \lin(\bold{B},\bold{C})$, where 
$\bold{B}$, $\bold{C}$ are formed by the columns of orthogonal matrices in the SVD-decomposition
of $\bold{S}$.
\end{proof}


In contrast to $\mathfrak{J}_{red}$ a minimizer of $\mathfrak{L}$ can 
be found rather easily. First, define 
\begin{equation}\label{centr}  
\bold{\bar{V}}^{C}:= \frac{\sum_{i=1}^{JK} \lambda_i \bold{\bar{V}}_i}{\sum_{i=1}^{JK} \lambda_i}
 \end{equation} 
as the {\bf centroid matrix }of $\bold{\bar{V}}_i$.

\begin{thm} \label{BCic}
Let $\bold{y_k}$, $\bold{z_k}$ be the left and right singular vectors in the SVD of $\bold{\bar{V}}^{C}$
defined in (\ref{centr}).
Then, 
\[ \bold{B_C} = (\bold{y_1} \ldots \bold{y_R}) \quad \bold{C_C} = (\bold{z_1} \ldots \bold{z_R})  \]
is a minimizer of $\mathfrak{L}(\bold{B},\bold{C})$.
Moreover, 
\[ \inf_{\bold{B},\bold{C}} \mathfrak{L}(\bold{B},\bold{C}) = \frac{1}{2} \left [   \left(1 - \|\bold{\bar{V}}^{C}\|_F^2 \right) \left ( \sum_{i=1}^{JK} \lambda_i \right )
 +   \left ( \sum_{i=1}^{JK} \lambda_i \right )
 \left ( \sum_{k = R+1}^{\min\{J,K\}} \sigma_k(\bold{\bar{V}}^{C})^2 \right) \right]
\]
\end{thm}
\begin{proof}
Expanding the square using $\|\bold{\bar{V}}_i\|_F^2 = 1$ yields
\begin{eqnarray*}
\sum_{i=1}^{JK} \lambda_i  \|\bold{\bar{V}}_i - \bold{X}\|_F^2 
&=& \left ( \sum_{i=1}^{JK} \lambda_i \|\bold{\bar{V}}_i\|_F^2\right )  - 
2 \left \langle \sum_{i=1}^{JK} \lambda_i \bold{\bar{V}}_i,\bold{X} \right \rangle  
+  \left (\sum_{i=1}^{JK} \lambda_i \right ) \langle \bold{X},\bold{X} \rangle\\  
&=& \left( \sum_{i=1}^{JK} \lambda_i \right) \left( 1- 2 \langle \bold{\bar{V}}^{C},\bold{X} \rangle 
+ \langle \bold{X},\bold{X} \rangle \right) \\
&=& 
 \left ( \sum_{i=1}^{JK} \lambda_i \right ) 
\left(1 - \|\bold{\bar{V}}^{C}\|_F^2  + \|\bold{\bar{V}}^{C}\|_F^2 - 
2 \langle \bold{\bar{V}}^{C}, \bold{X} \rangle + \langle \bold{X},\bold{X} \rangle \right )\\
&=&  \left( \sum_{i=1}^{JK} \lambda_i \right) \left(1 - \|\bold{\bar{V}}^{C}\|_F^2 \right) +
\left ( \sum_{i=1}^{JK} \lambda_i \right) \|\bold{\bar{V}}^{C}- \bold{X}\|^2.  
\end{eqnarray*}
Hence,
\[  \inf_{\bold{B},\bold{C}} \mathfrak{L}(\bold{B},\bold{C}) = \frac{1}{2} \left [
\left ( \sum_{i=1}^{JK} \lambda_i \right ) \left(1 - \|\bold{\bar{V}}^{C}\|_F^2 \right) + 
  \left ( \sum_{i=1}^{JK} \lambda_i \right )  \inf_{\mbox{rank}\{\bold{X}\}  \leq R} 
\|\bold{\bar{V}}^{C}- \bold{X}\|^2 \right ]
\]
Using again the Eckart-Young Theorem,  we see that a minimizer $\bold{X}$ is found through the truncated SVD of $\bold{\bar{V}}^{C}$,
\[ \bold{X}  = \sum_{k=1}^R \sigma_{k} \bold{y_k} \otimes \bold{z_k}. \]
Defining $\bold{B}$ and $\bold{C}$ as in the theorem we have that $X\in \lin(\bold{B},\bold{C})$
and thus ($\bold{B},\bold{C})$ yields a minimizer of $\mathfrak{L}(\bold{B},\bold{C})$.
\end{proof}
\begin{remark}
Computing minimizers of $\mathfrak{L}$ as in Theorem~\ref{BCic} 
yields matrices $\bold{B}$ and $\bold{C}$ which in turn approximate
 the minimizers of $\mathfrak{J}_{red}$. 
Theorem~\ref{main} yields also an simple algorithm using only linear algebra 
to calculate minimizers of $\mathfrak{L}$.
 We refer to this computation
of an initial guess as the {\bf Centroid Projection} algorithm. See Figure~\ref{fig:CentroidCP}a 
(Step~1--5)
for a detailed explanation of the implementation of the Centroid Projection algorithm.
\end{remark} 
Combining Corollary~\ref{clb} and Theorem~\ref{BCic} yields the following 
a-posteriori bounds on the quality of the output of the Centroid Projection algorithm. 

\begin{corollary}
Let $\bold{B_C}$ and $\bold{C_C}$ be as in Theorem~\ref{BCic}.
Then
\begin{equation}
 |\mathfrak{J}_{red}(\bold{B_C},\bold{C_C}) -
\inf_{(\bold{B},\bold{C})} \mathfrak{J}_{red}(\bold{B},\bold{C}) | \leq 
\frac{1}{2} \left( 
 \sum_{i=1}^{JK} \lambda_i  \left (\sum_{k=1}^{R} (\sigma_k^i)^2 \right )  
- \left( 
 \sum_{i=1}^{JK} \lambda_i \right) \left(
 \sum_{k = 1}^{R} \sigma_k(\bold{\bar{V}}^{C})^2 \right)
\right) 
\end{equation}
\end{corollary}
\begin{proof}
Note that the Frobenius norm can be expressed via the singular values
$\|\bold{\tilde{V}_i}\|_F^2 =  \sum_{k=1}^{JK} (\sigma_k^i)^2$. With the normalization
condition  $\|\bold{\tilde{V}_i}\|_F^2 = 1$,  Corollary~\ref{clb} and  Theorem~\ref{BCic}
the result follows.
\end{proof}
\begin{remark}
The possitivity of the right hand side in this estimate is a consequence of the
convexity of the sum of the square of the largest singular values (i.e. the Schatten norm).
\end{remark}
 
The output of the Centroid Projection algorithm can be used as 
sensible initial guesses for any current numerical methods for CP decomposition. 
Commonly, CP methods are initialized with random guesses which 
at times lead slow convergence rate. 
In Section~\ref{sec:num}, we describe how 
the Centroid Projection algorithm
is  able to mitigate the swamping effect which are often present
 in the ALS algorithm.  
which first 
\section{Numerical computation using the CPCP method}\label{sec:num}
The Centroid Projection algorithm yields an initial guess which in turn can
be combined with any iterative method for computing a CP approximation. 
We will for short refer to any combination of an iterative scheme using the
Centroid Projection an an initial guess as a {\bf CPCP} method. 
\subsection{CPCP with ALS schemes}\label{sec:numone}
Here we described some CP tensor decomposition numerical techniques based on the least-squares method.
 Matricizing $\mathcal{T} \approx \sum_{r=1}^R \bold{a_r} \circ \bold{b_r} \circ \bold{c_r}$ leads to 
three equivalent expressions:
\[\bold{T}^{JK \times I} \approx (\bold{B} \odot \bold{C}) \bold{A}^T,
 \bold{T}^{KI \times J} \approx (\bold{C} \odot \bold{A}) \bold{B}^T, 
\mbox{~and~} \bold{T}^{IJ \times K} \approx (\bold{A} \odot \bold{B}) \bold{C}^T.\]
To approximate the factors, three linear least-squares are solved iteratively:
\begin{eqnarray*}
\bold{A}^{k+1} &=& \min_{\bold{A}} \mathfrak{J}(\bold{A},\bold{B}^k,\bold{C}^k) \\ 
\bold{B}^{k+1} &=& \min_{\bold{B}} \mathfrak{J}(\bold{A}^{k+1},\bold{B},\bold{C}^k) \\ 
\bold{C}^{k+1} &=& \min_{\bold{C}} \mathfrak{J}(\bold{A}^{k+1},\bold{B}^{k+1},\bold{C}^k) \\ 
\end{eqnarray*}

\begin{itemize}
\item
{\bf Alternating Least-Squares (ALS) \cite{CarolChang,Harshman}.} 
$\mathfrak{J}(\bold{A},\bold{B}^k,\bold{C}^k)= \Vert \bold{T}^{JK \times I} - 
(\bold{B}^k \odot \bold{C}^k)  \Vert^2_F$, \\
$\mathfrak{J}(\bold{A}^{k+1},\bold{B},\bold{C}^k)$=$\Vert \bold{T}^{KI \times J} 
= (\bold{C} \odot \bold{A}) \bold{B}^T \Vert^2_F$ and 
$\mathfrak{J}(\bold{A}^{k+1},\bold{B}^{k+1},\bold{C}^k)$=$\Vert \bold{T}^{IJ \times K} -
 (\bold{A} \odot \bold{B}) \bold{C}^T \Vert^2_F$.
\item
{\bf Regularized Alternating Least-Squares (RALS) \cite{NaDeLatKind}.} 
$\mathfrak{J}(\bold{A},\bold{B}^k,\bold{C}^k)= \Vert \bold{T}^{JK \times I} - 
(\bold{B}^k \odot \bold{C}^k)  \Vert^2_F + 
\alpha_k \Vert  \bold{A}-\bold{A}^k\Vert^2_F$, $\mathfrak{J}(\bold{A}^{k+1},\bold{B},\bold{C}^k)$
=$\Vert \bold{T}^{KI \times J} - (\bold{C} \odot \bold{A}) \bold{B}^T \Vert^2_F + 
\alpha_k \Vert  \bold{B}-\bold{B}^k\Vert^2_F$ and\\
 $\mathfrak{J}(\bold{A}^{k+1},\bold{B}^{k+1},\bold{C}^k)$=
$\Vert \bold{T}^{IJ \times K} -(\bold{A} \odot \bold{B}) \bold{C}^T \Vert^2_F + 
\alpha_k \Vert  \bold{C}-\bold{C}^k\Vert^2_F$ where $\alpha_k$ is the regularization parameter.
\item
{\bf Rotationally Enhanced Alternating Least-Squares (REALS) \cite{PaaNaHop}.} 
The functional $\mathfrak{J}$ is similar to the ALS functional. However,
 \begin{eqnarray*}
\bold{A}^{k+1} \longleftarrow \bold{A}^{k+1} + d\bold{A}^k\\ 
\bold{B}^{k+1} \longleftarrow \bold{B}^{k+1}+ d\bold{B}^k\\ 
\bold{C}^{k+1} \longleftarrow \bold{C}^{k+1}+ d\bold{C}^k\\ 
\end{eqnarray*}
where  $d\bold{A}=\bold{A}\bold{R}$, $d\bold{B}=\bold{B}\bold{R}$, 
$d\bold{C}=\bold{C}\bold{R}$ and $\bold{R}$ is the rotational matrix.
\end{itemize}

\begin{figure}[htp]
  \begin{center}
   \subfigure[Centroid Projection Algorithm]{\label{fig:CentroidAlgo}
\includegraphics[width = 70 mm]{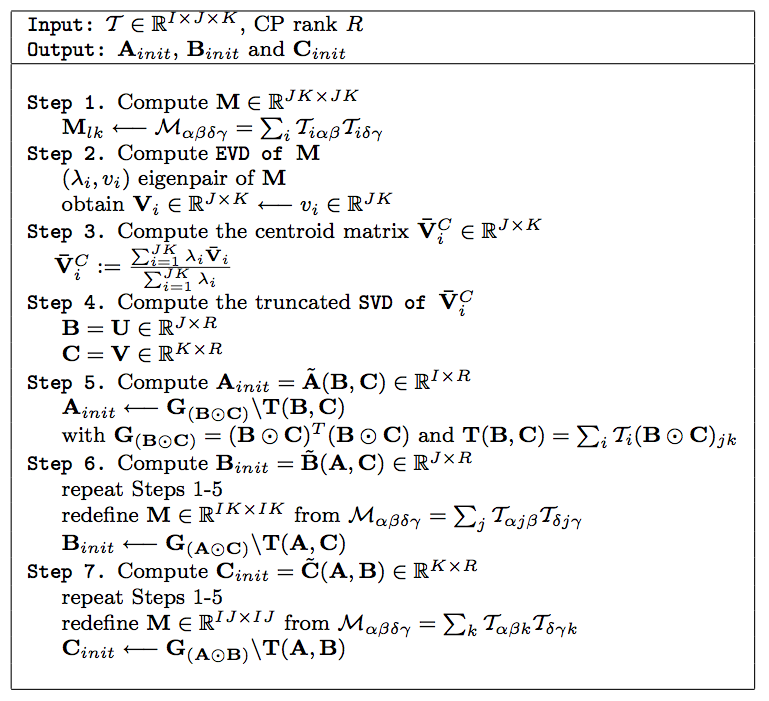}}
  \subfigure[Alternating CP Methods: (ALS, RALS, REALS) with Random Initial 
Conditions and (CALS, CRALS, CREALS) with Centroid Initial Conditions ]{\label{fig:CPMethods}
\includegraphics[width = 90 mm]{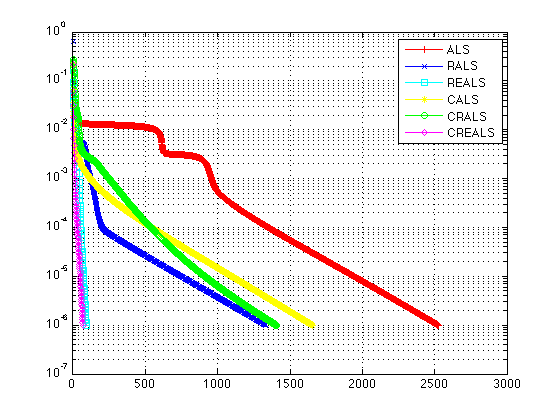}}
   \end{center}
  \caption{CP Methods with Random and Centroid Initial Conditions.}
  \label{fig:CentroidCP}
\end{figure}

The upper bound of $\mathfrak{J}_{red}$ in Theorem \ref{BCic} provides 
approximations for the factor matrices, closely estimating the solution subspaces. 
We called this method the Centroid Projection algorithm; it is summarized 
in Figure~\ref{fig:CentroidCP}a. Note that the Centroid Projection 
algorithm of Theorem~\ref{BCic} is contained in Step~1--5.
Steps~6--7 in Figure~\ref{fig:CentroidCP}a 
repeat the algorithm by interchanging the role of $\bold{A}$,
$\bold{B}$, $\bold{C}$. We observed a smaller initial residual error with this modification in most our numerical examples. The following CPCP methods, use initial conditions derived 
from the Centroid Projection and as an CP approximation 
the method ALS, RALS and REALS. We refer to them as 
 Centroid-ALS (CALS), Centroid-RALS (CRALS) 
and Centroid-REALS (CREALS), respectively. In Figure~\ref{fig:CentroidCP}b, 
we compared all six methods. Recall that a swamp is identified in a log error 
plot with an plateau and an extremely high number of iterations in order to converge. 
In most of our examples, both REALS and CREALS performed the fastest while
 ALS is the slowest, almost always hampered by a swamp. RALS, CRALS and CALS 
were comparable methods in performance, all dramatically decreasing the ALS swamp. 

The Centroid Projection method helps mitigate the effects of ALS swamps by providing a good set of initial factors lying close to the true solution subspaces.
\subsection{Symmetric CPCP}
\begin{figure}[htp]
  \begin{center}
    \subfigure[]{\label{fig:symCP1}\includegraphics[width = 70 mm]{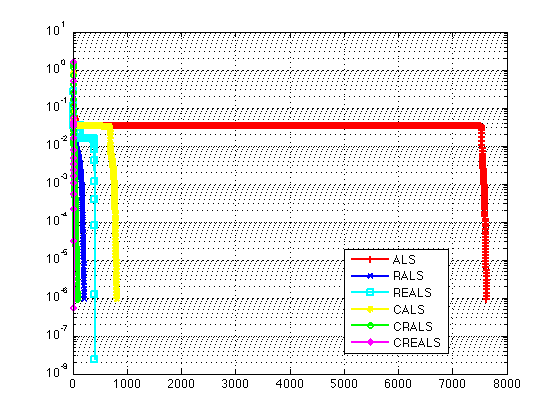}}
    \subfigure[]{\label{fig:symCP2}\includegraphics[width = 70 mm]{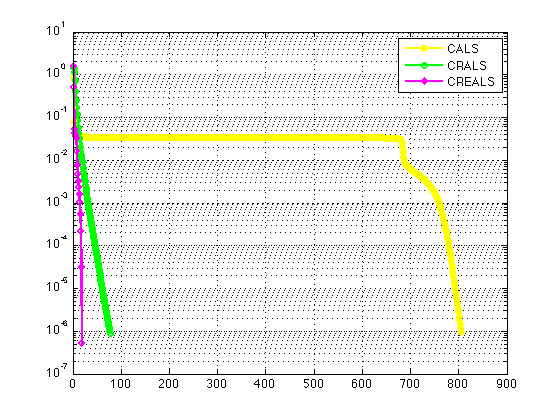}} 
\end{center}
  \caption{Fully Symmetric CP Decomposition.}
  \label{fig:SymmetricCP}
\end{figure}

The current methods for CP decomposition do not guarantee factorization with identical factors.  In fact, when ALS  is used in finding tensor decomposition with identical factors, the ALS algorithm will converge to a decomposition with no identical factors. Full and partial symmetries in tensor decomposition are referred to decomposition with at least two factors being identical. In a recent work of Stegeman \cite{Stegeman}, the existence and uniqueness of the $n$th order tensor decompositions with some form of symmetry have been studied for $n=3,4,5$. In other works \cite{BLNT,DeLatCasCar}, symmetries are also described by the permutation of the multi-indices of the tensor elements. An example is the following:  if $t_{ijk} = t_{jki}= t_{kij}$, then $\mathcal{T} \in \mathbb{R}^{N \times N \times N}$ with a tensor rank $R$ is fully symmetric and its factors are $\bold{A}=\bold{B}=\bold{C} \in \mathbb{R}^{N \times R}$. Another example is if $t_{ijkl}=t_{klij}$, then $\mathcal{T} \in \mathbb{R}^{M \times N \times M \times N}$ with a tensor rank $R$ is partially symmetric and its factors are $\bold{A}=\bold{C} \in \mathbb{R}^{M \times R}$ and $\bold{B}=\bold{D} \in \mathbb{R}^{N \times R}$.

When the Centroid Projection algorithm is applied to the CP methods (CALS, CRALS, CREALS) for symmetric decomposition, the methods with the centroid starters are guaranteed to converge to identical factors provided that the tensor dimensions and order satisfy the uniqueness and existence conditions of Kruskal \cite{Kruskal} and Stegeman \cite{Stegeman}. In the case that we have fully (partially) symmetric tensor, then $I=J=K$ ($I=J$ or $J=K$ or $I=K$). Thus,  from the EVD of $\bold{\bold{\bar{V}}_i} \in \mathbb{R}^{I \times I}$, we obtain the minimizers  $\bold{B}=\bold{C}$ in Theorem \ref{BCic} .

 Figure~\ref{fig:SymmetricCP} consists of plots of the number of iterations vs residual errors $\Vert \mathcal{T}_{orig}-\mathcal{T}_{est} \Vert_F^2$ for symmetric tensor decomposition with identical factors $\bold{A}=\bold{B}=\bold{C}$. For symmetric CP decomposition, CREALS has outperformed the other iterative methods 
in Section~\ref{sec:numone} up 
to a factor of $10^4$ while ALS has been consistently slow. 
ALS, RALS and REALS used random initial factors while CALS, CRALS 
and CREALS used calculated factors via the Centroid Projection algorithm. 
In most cases, these CPCP methods converged faster than the random-initialized CP methods.


\begin{thebibliography}{10}

%
%
%

\bibitem{BLNT}
M.~Brazell, N.~Li, C.~Navasca and C.~Tamon.
\newblock {\it Tensor and Matrix Inversions and Applications},
\newblock Submitted.


%
%
%
\bibitem{CarolChang}
J.D.~Carroll and J.J.~Chang. 
\newblock {\it Analysis of individual differences in multidimensional scaling via an N-way generalization of `Eckart-Young' decomposition,}
\newblock Psychometrika, {\bf 35} (1970),  pp.~283-319.
%
%
\bibitem{Comon}
P.~Comon, X.~Luciani and A.L.F.~de~Almeida.
\newblock {\em Tensor decompositions, altenating least squares and other tales.}
\newblock Journal of Chemometrics, {\bf 23} (2009) pp. 393-405.
%
%
%
%
%
%
\bibitem{HOOI}
L.~De Lathauwer, B.~De Moor and J.~Vandewalle.
\newblock {\it On the Best Rank-1 and Rank-(R1,R2,...,RN) Approximation of Higher-Order Tensors,}
\newblock SIMAX {\bf 21} 4 (2000), pp.~1324--1342.

\bibitem{DeLatCasCar}
L.~De Lathauwer, J.~Castaing, and J.-F.~Cardoso. 
\newblock {\it Fourth-Order Cumulant-Based Blind Identification of Underdetermined Mixtures.}
\newblock {IEEE Transactions on Signal Processing,} {\bf 55} 6 (2007), pp.~2965--2973.

\bibitem{SilLim}
V. de Silva and L.-H. Lim.
\newblock {\it Tensor rank and the ill-posedness of the best low-rank approximation problem.}
\newblock SIAM Journal on Matrix Analysis and Applications, {\bf 30} 3 (2008),  pp.~1084--1127.


%
%
%
%

\bibitem{Harshman}
R.A.~Harshman. 
\newblock {\it Foundations of the PARAFAC procedure: Models and conditions 
for an "explanatory" multi-modal factor analysis.}
\newblock {UCLA working papers in phonetics}, {\bf 16} (1970), pp.~1--84.

\bibitem{Ha00}
A.~Hatcher.
\newblock Algebraic topology.
\newblock Cambride Universtiy Press, Cambridge (2000).


\bibitem{Hitch1}
F.L.~Hitchcock. 
\newblock {\it The expression of a tensor or a polyadic as a sum of products.}
\newblock {Journal of Mathematics and Physics}, {\bf 6} (1927), pp.~164--189.

\bibitem{Hitch2}
F.L.~Hitchcock. 
\newblock {\it Multilple invariants and generalized rank of a p-way matrix or tensor.} 
\newblock {Journal of Mathematics and Physics}, {\bf 7} (1927), pp.~39--79.




%
%
%


%
%
%
%
%
\bibitem{Kolda}
T.~Kolda and B.W.~Bader.
\newblock {\it Tensor decompositions and applications,}
\newblock SIREV, {\bf 51} 3, (2009), pp.~455--500.

%
%

\bibitem{TheDutch}
W.P. Krijnen, T.K. Dijkstra and A. Stegeman.
\newblock {\it On the non-existence of optimal solutions and the occurrence of 
"degeneracy" in the Candecomp/Parafac model. }
\newblock Psychometrika, {\bf 73} (2008) pp.~431--439.



\bibitem{Kruskal}
J.B.~Kruskal.
\newblock {\it Three-way arrays: rank and uniquenss of 
trilinear decompositions with applications to arithmetic complexity and statistics,}
\newblock Linear Algebra and its Applications, {\bf 18} (1977), pp.~95--138.

%
%
%


\bibitem{LMV00}
L.~De~Lathauwer, B.~De~Moor, and J.~Vandewalle.
\newblock {\it On the best rank-1 and rank-$(R_1,R_2,. . .,R_N)$ approximation of
    higher-order tensors,}
\newblock SIAM J. Matrix Anal. Appl., {\bf 21} (2000), pp.~1324--1342.


\bibitem{LiKinNav}
N.~Li, S.~Kindermann and C.~Navasca.
\newblock {\it Some convergence results of a regularized alternating least-squares method for tensor decomposition,}
\newblock Submitted.


\bibitem{Mi63}
J.~Milnor.
\newblock {Morse theory.} 
\newblock Princeton University Press, Princton, 1963.




\bibitem{NaDeLatKind} 
C. Navasca, L. De Lathauwer and S. Kindermann.
\newblock {\em Swamp reducing technique for tensor decomposition,}
\newblock  in the 16th Proceedings of the European Signal Processing Conference, Lausanne, August 2008.
 
\bibitem{PaateroEngine}
P.~Paatero.
\newblock {\it The Multilinear Engine - a table-driven least squares program for solving
 multilinear problems, including the n-way Parallel Factor Analysis model.}
\newblock Journal of Computational and Graphical Statistics, {\bf 8} (1999), pp.~854--888.

 \bibitem{PP}
P.~Paatero.
\newblock Construction and analysis of degenerate PARAFAC models.
\newblock {\it  J. Chemometrics},  {\bf 14} (2000), pp.~285--299.  
 


\bibitem{PaaNaHop}
P.~Paatero, C.~Navasca, and P.~Hopke.
\newblock {\it Fast Rotationally Enhanced Alternating Least-Squares Method,}
\newblock Preprint. \texttt{http://people.clarkson.edu/$\sim$cnavasca/REALS.html}


%
%
%
%
%

\bibitem{RaCo}
M.~Rajih and P.~Comon. 
\newblock {\it Enhanced line search: A novel method to accelerate Parafac.}
\newblock in the 13th Proceedings of the European Signal Processing Conference, Antalya, Turkey, September 2005.


\bibitem{Stegeman}
A.~Stegeman.
\newblock {\it On uniqueness of the canonical tensor decomposition with some form of symmetry,} 
\newblock to appear in SIMAX.

\bibitem{StegDeLat}
A.~Stegeman and L.~De~Lathauwer.
\newblock  {\it A method to avoid diverging components in the Candecomp/Parafac model for generic
 $I\times J \times 2$ arrays.}
\newblock  SIAM Journal on Matrix Analysis and Applications, {\bf 30} (2009) pp.~1614--1638.
%
%
%
%
%
%
\bibitem{GR} G. Tomasi and R. Bro.
\newblock {\it A comparison of algorithms for fitting the PARAFAC model}
\newblock {Computational Statistics and Data Analysis,} {\bf 50} (2006), pp.~1700--1734.


%
%
%
%
%
%

\end{thebibliography}
\end{document}